\documentclass[preprint]{imsart}

\usepackage{amssymb}
\usepackage{amsfonts}
\usepackage{lineno}
\usepackage{graphicx}
\usepackage{caption}
\usepackage{subcaption}
\usepackage{float}
\usepackage{mathrsfs}
\usepackage{booktabs}

\startlocaldefs
\newtheorem{lemma}{Lemma}
\newtheorem{corollary}{Corollary}
\newtheorem{proposition}{Proposition}
\newtheorem{theorem}{Theorem}
\newtheorem{example}{Example}
\newtheorem{remark}{Remark}
\newcommand{\proof}{\noindent{\bf Proof.  }}
\newcommand{\eproof}{$\Box$}

\endlocaldefs

\begin{document}

\begin{frontmatter}

\title{Convergence Rates and Decoupling in Linear Stochastic Approximation Algorithms}
\runtitle{Convergence Rates and Decoupling}
 
\begin{aug}
\author{\fnms{Michael A.} \snm{Kouritzin}\ead[label=e1]{michaelk@ualberta.ca}},
\and
\author{\fnms{Samira} \snm{Sadeghi}\thanksref{t1}\ead[label=e2]{ssadeghi@ualberta.ca}}

\thankstext{t1}{Corresponding author. E-mail: ssadeghi@ualberta.ca}
\runauthor{M.A. Kouritzin and S. Sadeghi}

\affiliation{University of Alberta}

\address{Department of Mathematical and Statistical Sciences\\
University of Alberta, Edmonton, AB  T6G 2G1 Canada }

\end{aug}

\begin{abstract}
Almost sure convergence rates for linear algorithms 
$h_{k+1} = h_k +\frac{1}{k^\chi} (b_k-A_kh_k)$ are studied, 
where $\chi\in(0,1)$, $\{A_{k}\}_{k=1}^\infty$ are symmetric, positive semidefinite random matrices and 
$\{b_{k}\}_{k=1}^\infty$ are random vectors. 
It is shown that $|h_n- A^{-1}b|=o(n^{-\gamma})$ a.s. for the $\gamma\in[0,\chi)$, 
positive definite $A$ and vector $b$ 
such that $\frac{1}{n^{\chi-\gamma}}\sum\limits_{k=1}^n (A_{k}- A)\to 0$ 
and $\frac{1}{n^{\chi-\gamma}}\sum\limits_{k=1}^n (b_k-b)\to 0$ a.s. 
When $\chi-\gamma\in\left(\frac12,1\right)$, these assumptions are implied by the 
Marcinkiewicz strong law of large numbers, which allows the $\{A_k\}$ and $\{b_k\}$ to
have heavy-tails, long-range dependence or both.  
Finally, corroborating experimental outcomes and 
decreasing-gain design considerations are provided.
\end{abstract}



\end{frontmatter}
\section{Introduction}

Linear stochastic approximation algorithms have found widespread application in parameter
estimation, adaptive machine learning, 
signal processing, econometrics and pattern recognition 
(see, e.g., \cite{12}, \cite{13}, \cite{14}, \cite{15} and \cite{16}). 
Consequently, their asymptotic rates of almost sure and $r^{th}$-mean convergence as
well as invariance and large deviation principles are of utmost importance (see e.g., \cite{17}, 
\cite{18}, \cite{2}, \cite{19}, \cite{6}, \cite{20}, \cite{MK96}, \cite{11} and \cite{21}).
For motivation, suppose $\{x_k, k=1,2,\cdots\}$ and $\{y_k, k=2,3,\cdots\}$ are second order 
$\mathbb R^d-$ and $\mathbb R-$valued stochastic processes, defined on some probability 
space $(\Omega, F,P)$, that satisfy 
\begin{eqnarray}\label{observeqn}
y_{k+1}=x_k^Th+\epsilon_k,\quad\quad \forall k=1,2,\dots,
\end{eqnarray}
where $h$ is an unknown $d$-dimensional parameter or weight vector of interest and $\epsilon_k$ is a noise sequence.
One often wants to find the value of $h$ that minimizes the mean-square error
$h\to E|y_{k+1}-x_k^Th|^2$. 
This \emph{best} $h$ is given by $h=A^{-1}b$, 
where $A=E(x_kx_{k}^T)$ and $b=E (y_{k+1}x_k)$, assuming the expectations
exist, wide-sense stationarity conditions and that $A$ is positive definite.
However, we often do not know the joint distribution of $(x_k,y_{k+1})$ nor have the necessary
stationarity but instead estimate $h$ using a linear algorithm of the form: 
\begin{eqnarray}\label{a3}
	h_{k+1} = h_k +\mu_k (b_k-A_kh_k), 
\end{eqnarray}
where $\mu_k$ is the $k^{\rm th}$ step size (often of the form 
$\mu_k=k^{-\chi}$ for some $\chi\in\left(\frac12,1\right]$) and 
\begin{eqnarray}\label{a4}
A_k=\frac{1}{N} \sum_{l=\max\{k-N+1,1\}}^{k} x_lx_l^T,\,\mbox{and} \,\,
b_k=\frac{1}{N} \sum_{l=\max\{k-N+1,1\}}^{k}\! y_{l+1}x_l
\end{eqnarray}
for some $N\in\mathbb N $, are random sequences of symmetric, positive-semi-definite matrices
and vectors respectively.
Most often $N=1$ so $A_k=x_kx_k^T$ and $b_k=y_{k+1}x_k$.
More information on stochastic approximation can be found in 
e.g. \cite{27}, \cite{22}, \cite{5}, \cite{2}, \cite{3}, \cite{4} and \cite{1},
which provide examples and motivation for our work.
However, our work is easily differentiated from these.
Delyon \cite{27}, for example, focuses on non-linear stochastic approximation 
algorithms, treating linear examples the same as non-linear ones.
(In Section 4.2.2 he uses linear algorithm approximation but with a constant
deterministic matrix $A_k=A$ in our notation.)
Delyon's work handles important applications.
However, his A-stable and (A, B) Conditions 
are usually harder to verify than our Marcinkiewicz Strong Law of Large Numbers (MSLLN) conditions 
(given below) in the (unbounded, random $A_k$) linear case,
he does not supply almost sure rates of convergence, 
his theorems are geared to martingale-increment-plus-decreasing-perturbation noise
and he often assumes fourth order moments.
We are motivated by (but not restricted to) the common setting where $X_k^T=(x_k^T,y_{k+1})$ is a (multivariate) linear process
\begin{eqnarray}\label{xylinearproc}
X_k = \sum_{l=-\infty}^\infty C_{k-l}\Xi_l.
\end{eqnarray}
Matrix sequence $(C_l)$ can decay slowly enough (as $|l|\to\infty$) for long-range dependence (LRD)
while $\{\Xi_l\}$ can have heavy tails (HT), so $E|b_k|^2=\infty$ and/or $E|A_k|^2=\infty$.
Even in the lighter tail, short-range dependence case our two-sided linear process example $\{x_k\}$ is not a martingale.
Moreover; long-range dependence and heavy tails; exhibited in many network \cite{28}, financial and
paleoclimatic data sets for example; voids the usual mixing and moment conditions.
We focus on one-step versus Polyak-Ruppert's two-step averaging
algorithms but handle heavy tails and long range dependence, deriving a surprising
decoupling.
This means that the optimal convergence rate of (\ref{a3}) is affected by either
the heavy tails or the long-range dependence, whichever is worse, but not both. 
This contrasts the rate for partial sums of long-range dependent, heavily-tailed
random variables, which is degraded twice (see e.g. Theorem \ref{linear}).

Step size $\mu_k$ has a direct effect on the convergence rate and algorithm 
effectiveness (see, e.g \cite{23}, \cite{24} and references cited therein). 
Consider the extreme cases.
In the homogeneous, deterministic setting, i.e. $A_k=A$ and $b_k=b$, (\ref{a3}) can
solve the linear equation $Ah=b$ when matrix inversion of $A$ is ill-conditioned. 
In this case, a constant gain 
$\mu_k=\epsilon$ is best: 
Since $b=Ah$, we have $h_{k+1} = h_k-\epsilon A(h_k-h)$, so 
$h_{n}-h=(I-\epsilon A)^{n-1} (h_1-h)$ and 
$h_n\to h$ geometrically, provided $\epsilon$ is
small enough that the eigenvalues of $I-\epsilon A$ are within the unit disc. 
Conversely, in the presence of persistent noise, decreasing step sizes are required 
for the convergence $h_n\to h$. 
Existing results show that the best possible almost-sure rate of convergence 
is $|h_n-h|=O\left(\sqrt{n^{-1}\log\log(n)}\, \right)$, implied by the law of the iterated logarithm,
and that this rate is only attainable when $\mu_k=\frac1k$, 
second moments of $A_k,b_k$ exist and there is no long-range dependence. 
(These claims follow from the almost-sure invariance principle in Kouritzin \cite{MK96}.)

Herein, we handle all gains, long range dependence and heavy tails, addressing the optimal rate of convergence by establishing
results akin to the MSLLN, namely $|h_n-h|=o(n^{-\gamma})$ for all 
$\gamma<\gamma_0^{(\chi)}\doteq\chi-M$. $M$ is called the \emph{Marcinkiewicz threshold} in the sequel and is defined by
\small
\begin{eqnarray}\label{Mdef}
\quad M\doteq \inf\{\frac1m:
\lim_{n\to\infty} \frac{1}{n^{\frac1m}}\sum_{k=1}^{n}(A_k-A)=0 ,
\lim_{n\to\infty} \frac{1}{n^{\frac1m}}\sum_{k=1}^{n}(b_k-b)=0\ \mbox{ a.s.} \}. 
\end{eqnarray}
\normalsize
Usually, we expect $M\in (\frac12,1]$, 
due to Strong Law of Large Numbers and  Central Limit Theorem
 in the light-tail, short-range-dependence case but when there is LRD and/or HT $M$ generally cannot approach $\frac12$. 
When $\{(x_k^T,y_{k+1})^T:k\in\mathbb Z\}$ is a linear process as
in (\ref{xylinearproc}), it is shown in \cite{MS-small} that
$M= \frac1\alpha\vee\left(2-2\sigma\right)$
with $\alpha\doteq\sup\{a\le 2:\sup\limits_{t\ge 0}t^a P(|\Xi_1|^2>t)<\infty\}$ and
$\sigma\doteq\sup\{s\in(\frac12,1]:\sup\limits_l |l|^s \|C_l\|<\infty\}$. 
Hence, $\gamma<\gamma_0^{(\chi)}\doteq(\chi-\frac1\alpha)\wedge(\chi+2\sigma-2)$. 
Here, $\alpha\in(1,2]$ is a heavy-tail parameter with $\alpha=2$ indicating non-heavy tails
and $\sigma\in\left(\frac12,1\right]$ is a long-range dependence parameter with $\sigma=1$
indicating the minimal amount of long-range dependence.

In classical applications  
the best theoretical convergence rate is attained when $\chi=1$ corresponding to  $\gamma_0^{(\chi)}=\frac12$. 
However, this rate knowledge can lead to erroneous conclusions 
as the algorithm often performs better with
$\mu_k=k^{-\chi}$ for some $\chi<1$ than with $\mu_k=\frac1k$. 
How might one explain this apparent paradox?
First of all, these simple rate-of-convergence results do not account for the
possibility of exploding constants, i.e.\ if $h^{(\chi)}_k$ denotes the solution of
the algorithm (\ref{a3}) with $\mu_k=k^{-\chi}$, then 
$|h^{(\chi)}_n-h|=D^\chi n^{-\gamma^{(\chi)}}$ for all $\gamma^{(\chi)}<\gamma_0^{(\chi)}$.
However, this $D^\chi$ often increases rapidly as $\chi \nearrow 1$ so the \emph{observed}
convergence may be fastest for some $\chi<1$.
Secondly, a higher value of $\chi$ is worse for forgetting a poor initial guess $h_0$ of $h$
since you move further and further from the geometric convergence mentioned above as $\chi \to1$.

Our approach is to transfer the MSLLN from the partial sums of a linear algorithm's
coefficients to its solution.
In other words, we
establish the almost sure rates of convergence $|h_n-h|=o(n^{-\gamma})$ for the algorithm 
\begin{eqnarray}\label{a5}
h_{k+1} =h_k + \frac{1}{k^\chi} (b_k - A_kh_k)  \ \   \forall \ k=1,2,3,...
\end{eqnarray}
with $\chi\in(0,1)$, assuming only 
\begin{eqnarray}\label{a6}
\,\,\lim_{n\to\infty} \frac{1}{n^\chi}\sum_{k=1}^{n}(A_k-A)=0\ \mbox{ and }\,
\lim_{n\to\infty} \frac{1}{n^{\chi-\gamma}}\sum_{k=1}^{n}(b_k-A_kh)=0\ \mbox{ a.s.}
\end{eqnarray}
for some $\gamma \in [0,\chi)$, which can be implied by e.g.
\begin{eqnarray}\label{a6-1}
\lim_{n\to\infty} \frac{1}{n^{\chi-\gamma}}\sum_{k=1}^{n}(A_k-A)=0\ \mbox{ and }\
\lim_{n\to\infty} \frac{1}{n^{\chi-\gamma}}\sum_{k=1}^{n}(b_k-b)=0\ \mbox{ a.s.},
\end{eqnarray}
where $Ah=b$. 
When $\chi-\gamma \in (\frac{1}{2}, 1]$,
these conditions can be verified by the MSLLN under a variety of conditions,
which we study using the specific structure of $A_k$ and $b_k$ in Section 3.  

In addition to rates of convergence, our results show 
that convergence ($h_k\to h$) in (\ref{a5})
takes place provided that $\chi\in( M,1)$.  
All this suggests that more quickly decreasing gains like $\mu_k=\frac{1}{k^\chi}$ with $\chi$
near $1$ should be used in very heavy-tailed or long-range dependent settings.
Conversely, slowly deceasing gains like $\mu=\frac{1}{k^\chi}$ with smaller $\chi$ 
might work well in lighter-tailed, short-range-dependent situations.
Our simulations in Section 4 show that the smallest normalized error, 
$\frac{|h_{n}-h|}{|h_1-h|}$, usually occurs for $\chi \in (M ,1]$ 
and the most commonly used choice $\chi=1$ is most appropriate in very heavy-tailed
or long-range-dependent settings (where $M$ is close to $1$) or very long runs.
In other words, a slower decreasing gain usually gets you close to the true parameters $h$
more quickly unless the coefficients have a high probability of differing 
significantly from their means.

Let us consider what is new in terms of our theoretical results. 
The idea of inferring convergence and rates of convergence results for 
linear algorithms (\ref{a3})
from like convergence and rates of convergence of its coefficients is not new.
Indeed, it dates back at least to work done by one of the authors in 1992 and 1993 
(see \cite{6},\cite{20} and \cite{MK96}).
The first result \cite{6} considered relatively general gain $\mu_k$ and achieved
optimal rates of $r^{\rm th}$-mean convergence. 
It has been proved in \cite{MK96} that the solution of the linear algorithm (\ref{a3})
satisfies an almost sure invariance principle with respect to a limiting Gaussian process
when $\mu_k=\frac1k$ and each $A_k$ is symmetric
under the minimal condition that the coefficients satisfy such an a.s.\ invariance principle.
One could then immediately transfer functional laws of the iterated logarithm from
the limiting Gaussian process back to the solution of the linear algorithm. 
Again assuming the ``usual" conditions of $A_k$ symmetry and $\mu_k=\frac1k$, Kouritzin \cite{20}
showed that the solution of the linear algorithm converges almost surely given that the
coefficients do.
While this result does not state rates of convergence, our current work in going from
Proposition \ref{mainprop} to Theorem \ref{Converge} within shows that almost-sure rate of 
convergence sometimes follow from convergence results for linear algorithms as a simple 
corollary.

There were many results (see, e.g. \cite{5}, \cite{8} and \cite{7}) that preceded those mentioned above and
gave convergence or rates of convergence for linear algorithms.
However, these results assumed a specific dependency structure and, thereby,
were not generally applicable.
More recently, some authors, e.g. \cite{17}, \cite{27} and \cite{11}, have followed the path of 
transferring convergence and rates of convergence from partial sums
of (the coefficient) random variables to the solutions of linear equations.
Specifically, Tadi\'{c} \cite{11} transferred almost-sure rates of convergence, including
those of the law-of-the-iterated-logarithm rate, from the coefficients to the linear
algorithm in the non-symmetric-$A_k$, general-gain case. 
He does not develop a law of the iterated logarithm where one characterizes the
limit points nor does he consider functional versions.
Moreover, he imposes one of two sets of conditions (A and B in his notation).
Conditions B ensure the gain $\mu_k\approx\frac1k$, 
so these results should be compared to prior results in \cite{MK96} and 
\cite{B97}, which imply stronger Strassen-type functional laws of the iterated logarithm.
Tadi\'{c} does not give any examples verifying his Conditions A
where lessor rates are obtained.

It seems that  we are the first to consider processes that are simultaneously  heavy-tailed and long-range dependent in stochastic approximation.

The rest of this paper is organized as follows. Main theorems are formulated in Section $2$. Then,  
Section $3$ includes some background about the Marcinkiewicz Strong Law of Large Numbers for Partial Sums 
and a new MSLLN result for outer products of multivariate linear processes with LRD and HT. Experimental results 
are given in Section $4$ and proof of main result (Theorem \ref{Converge}) is delayed until Section $5$.

\section{Notation and Theoretical Result}
In this section, we define our notation and provide our results.
\subsection{Notation List}
$|x|$ is Euclidean distance of some $\mathbb R^d$-vector $x$.\\
$\|C\|=\sup_{|x|=1}|Cx|$ for any $\mathbb R^{n\times m}$-matrix $C$.\\ 
$\displaystyle\mid\mid\mid A\mid\mid\mid ^2 = \sum_{n=1}^{d} \sum_{o=1}^{d} (A^{(n,o)})^2$, 
$A^{(n,o)}$ is the $(n,o)^{\rm th}$ components of $A\in \mathbb R^{d\times d}$.\\ 
$\lfloor t\rfloor\doteq \max\{i\in\mathbb N_0:i\le t\}$ and
$\lceil t\rceil\doteq \min\{i\in\mathbb N_0:i\ge t\}$ for any $t\ge 0$.\\
$a_{i,k}\stackrel{i}\ll b_{i,k}$ means that for each $k$ there is a $c_k>0$ that does not
depend upon $i$ such that $|a_{i,k}|\le c_k |b_{i,k}|$ for all $i,k$.\\
$\prod\limits_{l=p}^q B_l$ ($\forall\, B_l$ being a $R^{d\times d}$-matrix) $=B_qB_{q-1}\cdots B_p$ if $q\ge p$ or $I$ if $p>q$.\\
$a\vee b=\max\{a,b\}$ and
$a\wedge b=\min\{a,b\}$.
\subsection{Main Results}
We will state and prove our results in a completely deterministic manner and then 
apply these results on a sample path by sample path basis. 
Therefore, we assume that $\chi \in (0,1)$, $d$ is a positive integer, $\{\bar{A}_{k}\}_{k=1}^\infty$ 
is a symmetric, positive semidefinite $R^{d\times d}$ -valued sequence, 
$\{\bar{b}_{k}\}_{k=1}^\infty$ is a $\mathbb R^d$-valued sequence and
$\{\bar{h}_k\}_{k=1}^\infty$ is a $\mathbb R^d$-valued sequence satisfying:
\begin{eqnarray}\label{b1}
  \bar{h}_{k+1} =\bar{h}_k + \frac{1}{k^\chi}(\bar{b}_k - \bar{A}_{k}\bar{h}_k ) \quad \mbox{for all } \quad k=1,2,...
\end{eqnarray}
Our first main result establishes rates of almost sure convergence:
\begin{theorem}\label{Converge}
Suppose $\gamma\in[0,\chi)$, $h\in \mathbb R^d$ and $A$ is a symmetric positive-definite matrix.\\
\noindent {\bf a)} If
\begin{eqnarray}\label{b2}
\lim_{n\to\infty} \left\|\frac{1}{n^\chi}\sum_{k=1}^n (\bar{A}_{k}- A)\right\|=0,\quad\mbox{and}
\end{eqnarray}
\begin{eqnarray}\label{cnver1}
\lim_{n \to \infty} 
\left|\frac{1}{n^{\chi-\gamma}} \sum_{k=1}^{n}( \bar{b}_{k}-\bar{A}_{k}h)\right|= 0,
\end{eqnarray}
then $|\bar{h}_n- h|=o(n^{-\gamma})$ as $n\to \infty$.\\

\noindent{\bf b)} Conversely, $\displaystyle \lim_{n \to \infty} \left|\frac{1}{n^\chi} 
\sum_{k=1}^{n} (\bar{b}_{k}-\bar{A}_{k}h)\right|= 0$, if $\lim\limits_{k\to\infty} 
\left|k^{1-\chi}(\bar{h}_k-h)\right|= 0$ and 
\begin{eqnarray}\label{b3}
\frac{1}{n^\chi} \sum_{k=1}^n k^{\chi-1}\left\|\bar{A}_{k}\right\|\quad \mbox{is bounded in n.}
\end{eqnarray}
\end{theorem}
\begin{remark} 
Lemma\ref{CondComp} of Appendix establishes that (\ref{b2}) implies (\ref{b3}). 
\end{remark}
Theorem \ref{Converge} with $\bar{A}_{k}=A_k(\omega)$ and $\bar{b}_k 
= b_k(\omega)$ for all $k$, implies $h_n(\omega)$, the solution of (\ref{a3}), 
converges to $h = A^{-1}b$ a.s.
Indeed, to establish the rate of convergence $|\bar{h}_n- h|=o(n^{-\gamma})$, one need only
check standard conditions for the MSLLN in (\ref{b2}) and (\ref{cnver1}),
which is less onerous task than checking the technical conditions 
in Corollary 1 or Corollary 3
in \cite{11} say.
Indeed, there appears to be a need for some extra stability in \cite{11} by the imposition that 
``the real parts of the eigenvalues of $A$ should be strictly less than a certain negative value
depending on the asymptotic properties of $\{\gamma_n\}$ and $\{\delta_n\}$".
We do not need any such extra condition. 

Generally, we do not know $h$ when using stochastic approximation so we cannot just verify (\ref{cnver1})
but rather use the following corollary instead of Theorem \ref{Converge}.
\begin{corollary}\label{CorCon}
Suppose $\gamma\in[0,\chi)$ and $A$ is a symmetric positive-definite matrix.
\begin{eqnarray}\label{cnver}
\frac{1}{n^{\chi-\gamma}} \sum_{k=1}^{n}( \bar{b}_{k}-b)\to 0 \mbox{ and } 
\frac{1}{n^{\chi-\gamma}} \sum_{k=1}^{n}( \bar{A}_{k}-A)\to 0\ \mbox{a.s.}
\end{eqnarray}
Then, $|\bar{h}_n- h|=o(n^{-\gamma})$ as $n\to \infty$. 
\end{corollary}

Finally, we give a version of the theorem for linear processes under very verifiable conditions.
\begin{theorem}\label{LinearMain}
Let $ \left\{\Xi_{l}\right\}$ be i.i.d.\ zero-mean random $\mathbb R^m$-vectors such that 
$$\sup\limits_{t\ge 0 }t^{\alpha}P(|\Xi_1|^2>t)<\infty\quad \mbox{for some}\,\, \alpha\in(1,2)$$
$(C_{l})_{l\in\mathbb Z}$ be $\mathbb R^{(d+1)\times m}$-matrices such that
$\sup\limits_{l\in\mathbb Z}|l|^\sigma\|C_l\|<\infty\, \mbox{for some}\, \sigma\in \left(\frac{1}{2},1\right]$,    
\begin{eqnarray}\nonumber
(x_k^T,y_{k+1})^T = \sum_{l=-\infty}^\infty C_{k-l}\Xi_l,
\end{eqnarray}
$A_k=x_kx_k^T$, $b_k=y_{k+1}x_k$ and
$A=E[x_kx^T_k]$ and $b=E[y_{k+1}x_k]$.\\
Then, $|\bar{h}_n- h|=o(n^{-\gamma})$ as $n\to \infty$ a.s. for any 
$\gamma<\gamma_0^{(\chi)}\doteq(\chi-\frac1\alpha)\wedge(\chi+2\sigma-2)$.
\end{theorem} 
\begin{remark}
Theorem \ref{LinearMain} follows from Corollary \ref{CorCon} and Theorem \ref{MatrixHeavyLong} (to follow), by letting
$\frac1p = \chi-\gamma$ and $\overline{X}_k^T=X_k^T=(x_k^T,y_{k+1})$ and correspondingly, $\overline{\Xi}_{l}=\Xi_l$, 
$\overline{C}_{l}=C_{l}$ and $\overline{\sigma}=\sigma$. $\sigma$ and $\alpha$ are long-range dependence and heavy-tail parameters, respectively.
Theorem \ref{MatrixHeavyLong} also appears in \cite[Theorem $4$]{MS-small}.
\end{remark}

\section{Marcinkiewicz Strong Law of Large Numbers for Partial Sums}
Our basic assumptions are MSLLN for random variables for $\{A_k\}$ and $\{b_k\}$.
(Technically, our assumptions are even more general as they allow the non-MSLLN case where
$\chi-\gamma\le\frac12$ that could be verified by some other method in some
special situations.)
The beauty of this MSLLN assumption is that: 1) It is minimal in the sense
that the linear algorithm with $A_k=I$ and $\mu_k=\frac1k$ reduces to
the partial sums
\(
h_{k+1}-h =\frac1k \sum\limits_{j=1}^k (b_j-b) 
\)
(since $h=b$ when $A=I$) so a rate of convergence in the algorithm solution $h_k$ implies
a MSLLN for random variables $\{b_j\}$.
2) MSLLNs hold under very general conditions, including heavily-tailed
and long-range dependent data.
Hence, we review some of the literature in this area before giving simulation results for our theoretical work.

The classical independent case, due to Marcinkiewicz, is 
generalized slightly by Rio \cite{29}:
\begin{theorem}\label{mdependent}
Let $\{X_i\}$ be an $m$-dependent, identically distributed sequence of 
zero-mean $\mathbb R$-valued random variables such that $E|X_1|^p<\infty$ for some
$p\in(1,2)$. 
Then, 
\[
\frac{1}{n^\frac1p}\sum_{i=1}^n X_i\to 0\ \mbox{a.s.}
\]
\end{theorem}
Actually, Rio gives a more general $m$-dependent result on page 922 of his work.
However, the important observation for us is that 
only the $p^{\rm th}$ moment need be finite rather than a higher
moment as is typical under some stronger dependence assumptions.
Theorem \ref{mdependent} is quite useful in verifying our conditions 
when $\{A_k\}$ and $\{b_k\}$ may
have heavily-tailed distributions but are independent or $m$-dependent.
For example, if $\chi-\gamma\in\left(\frac12,1\right)$ and
the $\{A_k\}$ and $\{b_k\}$ are defined as in (\ref{a4}) in terms of
i.i.d. $\{x_k\}$ and $\{y_k\}$ with $E|x_1|^{\frac{2}{\chi-\gamma}}<\infty$,
$E|y_1|^{\frac{2}{\chi-\gamma}}<\infty$, then $\{A_k\}_{k\ge M}$ and $\{b_k\}_{k\ge M}$ are 
identically distributed, $M$-dependent and
\[
\frac1{n^{\chi-\gamma}}\sum_{k=1}^n(A_k-A)\to 0\ \mbox{ and }\
\frac1{n^{\chi-\gamma}}\sum_{k=1}^n(b_k-b)\to 0
\]
a.s., where $A=EA_k$ and $b=Eb_k$, by applying Theorem \ref{mdependent} for each
component. 
Hence, (\ref{a6-1}) holds.

There are many other important results that include heavy-tails, long-range 
dependence or both.
For example, Louhchi and Soulier \cite{30} prove the following result for
linear symmetric $\alpha$-stable (S$\alpha$S) processes. 
\begin{theorem}\label{linear}
Let $\{\zeta_j\}_{j\in\mathbb Z}$ be i.i.d. sequence of S$\alpha$S random variables
with $1<\alpha<2$ and $\{c_j\}_{j\in\mathbb Z}$ be a bounded collection such that 
\(
\sum\limits_{j\in\mathbb Z}|c_j|^s<\infty
\)
for some $s\in[1,\alpha)$. 
Set
\(
X_k=\sum\limits_{j\in\mathbb Z}c_{k-j}\zeta_j.
\)
Then, 
for $p\in(1,2)$ satisfying $\frac1p>1-\frac1s+\frac1\alpha$ 
\[
\frac1{n^\frac1p}\sum_{i=1}^n X_i\to 0\ \mbox{a.s.}
\]
\end{theorem}
The condition $s<\alpha$ ensures 
\(
\sum\limits_{j\in\mathbb Z}|c_j|^\alpha<\infty
\)
and thereby convergence of 
\(
\sum\limits_{j\in\mathbb Z}c_{k-j}\zeta_j.
\)
Moreover, $\{X_k\}$ not only exhibits heavy tails but also long-range dependence
if, for example, $c_j=|j|^{-\sigma}$ for $j\neq 0$ and some $\sigma\in \left(\frac12,1\right)$.
Notice there is interactions between the heavy tail condition and the long range dependent 
condition.
In particular for a given $p$, heavier tails ($\alpha$ becomes smaller) 
implies that you cannot have as 
long range dependence ($s$ becomes smaller) and vice versa.
Moreover, this result is difficult to apply in the stochastic approximation 
setting.
For example, if wanted to apply it for $X_k=A_k$ in the scalar case, then
we would need $x_k$ such that $x_k^2=A_k$ which is impossible when $A_k$
is S$\alpha$S.

One nice feature of mixing assumptions is that they usually transfer from
random variables to functions (like squares) of random variables.
There are many mixing results that handle long range dependence.
For example, Berbee \cite{31} gives a nice $\beta$-mixing result.
However, strong mixing is one of the most general types of mixing that
is more easily verified in practice.
Hence, we will just quote the following strong mixing result from Rio \cite{29} (Theorem 1)
in terms of the inverse $\alpha^{-1}(u)=\sup\{t\in\mathbb R^+: \alpha_{\lfloor t\rfloor}>u\}$
of the strong mixing coefficients
\[
\alpha_n=\sup_{k\in\mathbb Z}
\sup_{A\in\sigma(X_i, i\le k-n),B\in\sigma(X_k)}|P(AB)-P(A)P(B)|
\]
and the complementary quantile function
\[
Q_X(u)=\sup\{t\in\mathbb R^+: P(|X|>t)>u\}. 
\]
\begin{theorem}\label{strongmixing}
Let $\{X_i\}$ be an identically-distributed zero-mean sequence of 
$\mathbb R$-valued random variables such that 
\(
\int_0^1 [\alpha^{-1}(t/2)]^{p-1}Q_X^p(t)dt<\infty
\)
for some $p\in(1,2)$. 
Then, 
\[
\frac1{n^\frac1p}\sum_{i=1}^n X_i\to 0\ \mbox{a.s.}
\]
\end{theorem}
Notice again that for a given $p$, heavier tails implies that you cannot have as 
long range dependence and vice versa:
If you wanted to maintain the same value of the integral condition and
there became more area under $P(|X|>t)$, then there would be 
more area under $Q_X^p(t)$ so the area under $[\alpha^{-1}(t/2)]^{p-1}$, which
is equal to $2\sum\limits_{n=0}^\infty \alpha_n^{p-1}$, would have to decrease
to compensate. Also, there can be difficulty in establishing that a given model
satisfies the strong mixing condition with the required decay of mixing coefficients.
Still, this is an important result for verifying our basic assumptions.

A new MSLLN for outer products of multivariate linear processes with long range dependence and heavy tails 
is studied in \cite{MS-small}. A new decoupling property is proved that shows the convergence rate 
is determined by the worst of the heavy tails
or the long range dependence, but not the combination.
This result used to obtain Marcinkiewicz Strong Law of Large Numbers
for stochastic approximation (Theorem \ref{LinearMain}). The result is as follow.
\begin{theorem}\label{MatrixHeavyLong}
Let $ \left\{\Xi_{l}\right\}$ and $ \left\{\overline{\Xi}_{l}\right\}$ be i.i.d.\ zero mean random $\mathbb R^m$-vectors such that 
$\Xi_{l}=\left( \xi_l^{(1)}, . . ., \xi_l^{(m)}\right) $, $\overline{\Xi}_l=\left( \overline{\xi}_l^{(1)}, . . ., \overline{\xi}_l^{(m)}\right) $, 
$E[|\Xi_1|^2]<\infty$, $E[|\overline{\Xi}_1|^2]<\infty$ and
 $\max\limits_{1\le i,j\le m}\sup\limits_{t\ge 0}t^{\alpha}P(|\xi_1^{(i)}\overline{\xi}^{(j)}_1|>t)<\infty$
for some $\alpha \in(1,2)$.
Moreover, suppose matrix sequences
$(C_{l})_{l\in\mathbb Z}, (\overline{C}_{l})_{l\in\mathbb Z}$ $\in \mathbb R^{(d+1)\times m}$ satisfy
\begin{eqnarray*}
\sup\limits_{l\in\mathbb Z}|l|^\sigma\|C_l\|<\infty,\, \sup\limits_{l\in\mathbb Z}|l|^{\overline{\sigma}}\|\overline{C}_l\|<\infty \quad 
\mbox{for some} \quad (\sigma,\overline{\sigma}) \in \left(\frac{1}{2},1\right], 
\end{eqnarray*}
$X_k$, $\overline{X}_k$ take form of (\ref{xylinearproc}), $D_k=X_k\overline{X}_k^T$ and $D=E[X_1\overline{X}_1^T]$.
Then, for $p$ satisfying $p<\frac{1}{2-\sigma -\overline \sigma }\wedge\alpha$ 
\[
\lim\limits_{n\rightarrow \infty}\frac{1}{n^{\frac{1}{p}}}\sum\limits_{k=1}^{n}\left(D_{k}-D\right)=0\ \quad \mbox{a.s.} 
\]
\end{theorem}
 
This theorem actually shows the MSLLN for $D_k-E[D_k]$, where
\(
D_k=\left(\begin{array}{l c r}
	x_kx_k^T&&y_{k+1}x_k\\
	y_{k+1}x_k^T&&y_{k+1}^2
\end{array}\right)
\),
and then throw out the unneeded columns.
\section{Experimental Results}  
In this section we now verify our results of the previous section experimentally in the stochastic 
approximation setting discussed in the introduction.
In particular, we use
\emph{power law} or \emph{folded t} distributions.

{\bf Power law distribution:} A random variable $\xi$ obeys a power law with parameters 
$\beta>1$ and $x_{\min}>0$, written $\xi\sim PL (x_{min},\beta)$, if it has density 
\begin{eqnarray*}
f(x)=\frac{\beta-1}{x_{min}} (\frac{x}{x_{min}})^{-\beta} \quad \forall \,\, x\geq x_{min}
\end{eqnarray*}
Note that 
$E|\xi|^r= \left\{\begin{array}{ll} x_{min}^r (\frac{\beta-1}{\beta-1-r})& \quad r<\beta-1\\
\infty & \quad r\ge\beta-1\end{array}\right.$.

{\bf Folded t distribution:} A non-negative random variable $\xi$ has a folded $t$ distribution with 
parameter $\beta>1$, written $\xi\sim Ft (\beta)$, if it has density  
\begin{eqnarray*}
f(x)=\frac{2\Gamma( \frac{\beta}{2})}{\Gamma( \frac{\beta-1}{2})\sqrt{(\beta-1) \pi}} \left( 1+\frac{x^2}{(\beta-1)} \right) ^{-\frac{\beta}{2}} \quad \forall \,\, x> 0.
\end{eqnarray*}
Note that $E(|\xi|^r)$ exists if and only if $r<\beta-1$. 

Experimental results in this section are divided in two parts.  
\subsection{Heavy-tailed cases}
Assume $N=1$ in (\ref{a4}), dimension is $d=2$ and 
$\{(x_k^{(1)},x_k^{(2)},\epsilon_k)^T,\ k=1,2,...\}$ are i.i.d.\ random vectors 
so linear algorithm (\ref{a3}) reduces to:
	\begin{equation}\label{a-2}
	\quad \quad h_{k+1} = h_k +\mu_k (x_k y_{k+1}-x_kx_{k}^Th_k)=h_k +\mu_k (x_k x_{k}^Th+x_k\epsilon_k-x_kx_{k}^Th_k).		
	\end{equation}
For consistency and performance, we always let $x_k^{(1)},x_k^{(2)}$ and $\epsilon_k$ be independent. 
The runs are always initialized with $h_1=(101,101)^T$ and, for testing purposes, 
the optimal $h=(1,1)^T$ is known. 
\begin{example}
Let $x_k^{(1)},x_k^{(2)} \sim PL(x_{min}=1,\beta)$ and 
$\epsilon_k=\epsilon'_k- E(\epsilon'_k)$ with $\epsilon'_{k} \sim PL(x'_{min}=0.01,\beta)$. 
The normalized errors in 100 trial simulations, $\{h_n^{(i)}\}_{i=1}^{100}$, are averaged 
\(
\displaystyle \overline{\mbox{rh}}= \frac{1}{100}\sum_{i=1}^{100} \frac{|h_{n}^{(i)}-h|}{|h_1-h|}
\)
and given in the Table~\ref{tab:b-1} in terms of gain parameter $\chi$,
distributional parameter $\beta$ and sample size $n$. 
 
\begin{table}[h!]
\caption{Algorithm performance-Power Law} 
\centering 
\small\addtolength{\tabcolsep}{-2pt}
\begin{tabular}{l c c c c c c c c c c c}  
\hline\hline      
&\multicolumn{3}{c}{ n=100000} & & \multicolumn{3}{c}{ n=750000} & &\multicolumn{3}{c}{ n=1500000} \\[-0.15ex]
\cline{3-3} \cline{7-7} \cline{11-11}              
$\chi \backslash \beta$ & 3.5  & 4   & 4.5 &  & 3.5   & 4   & 4.5  & & 3.5  & 4   & 4.5 \\ [-0.15ex]
\cline{2-4} \cline{6-8} \cline{10-12}
0.55  & 0.1043 & 0.0391  & 0.0214  &   & 0.0841  & 0.0315  & 0.0155  &   & 0.0625  & 0.0254   & 0.0134  \\ [-0.15ex]
0.6   & 0.0864 & 0.0314  & 0.0169  &   & 0.0707  & 0.0243  & 0.0115  &   & 0.0548  & 0.0203   & 0.0099  \\ [-0.15ex]
0.65  & 0.0690 & 0.0247  & 0.0129  &   & 0.0578  & 0.0192  & 0.0086  &   & 0.0469  & 0.0166   & 0.0075  \\ [-0.15ex]
0.7   & 0.0525 & 0.0190  & 0.0098  &   & 0.0487  & 0.0159  & 0.0067  &   & 0.0457  & 0.0141   & 0.0056  \\ [-0.15ex]
0.75  & 0.0397 & 0.0151  & 0.0082  &   & 0.0449  & 0.0137  & 0.0051  &   & 0.0456  & 0.0114   & 0.0042  \\ [-0.15ex]
0.8   & 0.0326 & 0.0136  & 0.0105  &   & 0.0448  & 0.0111  & 0.0038  &   & 0.0402  & 0.0087   & 0.0031  \\[-0.15ex]
0.85  & 0.0314 & 0.0168  & 0.0549  &   & 0.0398  & 0.0085  & 0.0082  &   & 0.0324  & 0.0070   & 0.0035  \\[-0.15ex]
0.9   & 0.0344 & 0.0719  & 0.2445  &   & 0.0438  & 0.0118  & 0.0764  &   & 0.0272  & 0.0079   & 0.0341  \\[-0.15ex]
0.95  & 0.0902 & 0.3047  & 0.6631  &   & 0.3739  & 0.0897  & 0.3068  &   & 0.0248  & 0.0519   & 0.1963  \\[-0.15ex]
0.98  & 0.2226 & 0.5733  & 1.0154  &   & 0.9219  & 0.2302  & 0.5251  &   & 0.0374  & 0.1488   & 0.3930  \\[-0.15ex]
1     & 0.3876 & 0.8062  & 0.6631  &   & 1.1891  & 0.3745  & 0.6925  &   & 0.0662  & 0.2596   & 0.5644  \\[-0.15ex]
\hline\hline 
\end{tabular}
\label{tab:b-1}
\end{table}

The Marcinkiewicz threshold, $M=\frac{2}{(\beta-1)}$, corresponding to $\beta=3.5$, $\beta=4$ and $\beta=4.5$ are
respectively $M=0.8$, $0.67$ and $0.57$.
Our theoretical results prove convergence above this threshold.
While the results in Table 1 are obviously still influenced by (heavy-tailed) randomness, one can see that
convergence does appear to be taking place as one moves from $n=100,000$
through $n=750,000$ to $n=1,500,000$ when $\chi>M$ and it is less clear
that convergence is taking place when $\chi<M$.
Furthermore, our (as well as prior) theoretical
results predict rates of convergence that increase in $\chi$.
Indeed, in the case $\beta=4$ our theoretical results suggest that 
$\chi\approx 1$ should result in a rate of convergence $|h_n-h|=o(n^{-0.33})$
while $\chi= 0.85$ should only result in a rate of convergence $|h_n-h|=o(n^{-0.18})$.
Conversely, Table ~\ref{tab:b-1} demonstrates that $\chi=0.85$ performs better, which seems to contradict the theory.
However, this paradox is explained by the exploding constants discussion of the introduction
and, in fact, points out that more refined theory, involving functional results, is needed.
The proper way to use our theoretical results then is to predict 
the best $\chi$ (lowest value of $\overline{rh}$) in the range of $(M,1]$ i.e.\ in
$(0.8,1]$, $(0.67,1]$ and $(0.57,1]$, respectively for our three $\beta$'s.

\begin{table}[h!]
\caption{Best fixed $\chi$-Power Law} 
\centering  
\begin{tabular}{c c c c c c c c c c c c}  
\hline\hline      
&\multicolumn{3}{c}{ n=100000} & & \multicolumn{3}{c}{ n=750000} & &\multicolumn{3}{c}{ n=1500000} \\   [-0.5ex]
\cline{3-3} \cline{7-7} \cline{11-11}     
$ \beta$ & 3.5  & 4   & 4.5 &  & 3.5   & 4   & 4.5  & & 3.5  & 4   & 4.5 \\   [-0.5ex]
\cline{2-4} \cline{6-8} \cline{10-12}
 \small{Best} $\chi$ &0.85 &0.8  &0.75& &0.85 &0.85 &0.8 & &0.95 &0.85 &0.8 \\ [-0.5ex]
 \small{Resulting} $\gamma$              &0.05 &0.13 &0.18& &0.05 &0.18 &0.23& &0.15 &0.18 &0.23 \\ [-0.5ex]
\hline\hline 
\end{tabular}
\label{tab:b-2}
\end{table}

The best $\chi$'s, corresponding to the smallest value of $\overline{\mbox{rh}}$ for 
$\beta=3.5$, $\beta=4$ and $\beta=4.5$ and 3 different sample sizes, 
as well as the $\gamma$ corresponding to the theoretical rate of convergence $o(n^{-\gamma})$
are summarized in Table ~\ref{tab:b-2}.
In all cases the best value for $\chi$ is in the predicted range. 
As we explained, a faster decreasing gain is appropriate for a heavier-tailed distribution,
which is also confirmed by Table 2.
Notice also that the best $\chi$ increases in $n$, a phenomenon consistent with
our exploding constants and the initial condition effect discussion. 
\end{example}
Now, we repeat the previous example with a different distribution.
Since the results are consistent with those of the previous example, we will
keep our discussion to a minimum.
\begin{example}
Let $x_k^{(1)},x_k^{(2)} \sim Ft(\beta)$ and 
$\epsilon_k=\epsilon'_k- E(\epsilon'_k)$ with $\epsilon'_{k} \sim Ft(\beta)$. 
The simulation results for three $\beta$'s: $3.5, 4$ and $4.5$ with corresponding 
Marcinkiewicz thresholds, $M=\frac2{\beta-1}$, $0.8$, $0.67$ and $0.57$ are given in Table ~\ref{tab:tb-1} with 
sample sizes: $n=50,000, 100,000$ and $750,000$.  
\begin{table}[h!]
\caption{Algorithm performance-Folded t}
\centering  
\small\addtolength{\tabcolsep}{-2pt}
\begin{tabular}{l c c c c c c c c c c c }  
\hline\hline      
&\multicolumn{3}{c}{ n=50000} & & \multicolumn{3}{c}{ n=100000} & &\multicolumn{3}{c}{ n=750000} \\   [-0.15ex]  
\cline{3-3} \cline{7-7} \cline{11-11}         
$\chi \backslash \beta$ & 3.5  & 4   & 4.5 &  & 3.5   & 4   & 4.5  & & 3.5  & 4   & 4.5 \\   [-0.15ex]
\cline{2-4} \cline{6-8} \cline{10-12}
0.55  & 0.1040 & 0.0422  & 0.0274  &   & 0.1003  & 0.0429  & 0.0221  &   & 0.0780  & 0.0246   & 0.0145  \\ [-0.15ex]
0.6   & 0.0958 & 0.0345  & 0.0221  &   & 0.0929  & 0.0336  & 0.0177  &   & 0.0590  & 0.0195   & 0.0104  \\[-0.15ex]
0.65  & 0.0851 & 0.0291  & 0.0177  &   & 0.0778  & 0.0269  & 0.0141  &   & 0.0420  & 0.0149   & 0.0081  \\[-0.15ex]
0.7   & 0.0697 & 0.0245  & 0.0138  &   & 0.0661  & 0.0216  & 0.0112  &   & 0.0318  & 0.0120   & 0.0064  \\[-0.15ex]
0.75  & 0.0599 & 0.0204  & 0.0113  &   & 0.0556  & 0.0173  & 0.0089  &   & 0.0336  & 0.0099   & 0.0050  \\[-0.15ex]
0.8   & 0.0505 & 0.0172  & 0.0103  &   & 0.0439  & 0.0140  & 0.0075  &   & 0.0374  & 0.0076   & 0.0038  \\[-0.15ex]
0.85  & 0.0399 & 0.0145  & 0.0098  &   & 0.0341  & 0.0118  & 0.0063  &   & 0.0339  & 0.0058   & 0.0029  \\[-0.15ex]
0.9   & 0.0312 & 0.0133  & 0.0087  &   & 0.0278  & 0.0100  & 0.0057  &   & 0.0265  & 0.0048   & 0.0024  \\[-0.15ex]
0.95  & 0.0275 & 0.0241  & 0.0097  &   & 0.0245  & 0.0089  & 0.0060  &   & 0.0205  & 0.0039   & 0.0021  \\[-0.15ex]
0.98  & 0.0347 & 0.0475  & 0.0212  &   & 0.0274  & 0.0117  & 0.0121  &   & 0.0179  &0.00371   & 0.0032  \\[-0.15ex]
0.99  & 0.0404 & 0.0583  & 0.0295  &   & 0.0310  & 0.0149  & 0.0172  &   & 0.0173  &0.00373   & 0.0048  \\[-0.15ex]
1     & 0.0486 & 0.0700  & 0.0413  &   & 0.0369  & 0.0205  & 0.0249  &   & 0.0170  & 0.0039   & 0.0077  \\[-0.15ex]
\hline\hline 
\end{tabular}
\label{tab:tb-1}
\end{table}

\begin{table}[h!]
\begin{center}  
\caption{Best fixed $\chi$-Folded t}
\begin{tabular}{c c c c c c c c c c c c}  
\hline\hline      
&\multicolumn{3}{c}{ n=50000} & & \multicolumn{3}{c}{ n=100000} & &\multicolumn{3}{c}{ n=750000} \\   [-0.5ex]
\cline{3-3} \cline{7-7} \cline{11-11}    
$ \beta$ & 3.5  & 4   & 4.5 &  & 3.5   & 4   & 4.5  & & 3.5  & 4   & 4.5 \\   [-0.5ex] 
\cline{2-4} \cline{6-8} \cline{10-12}
\small{Best} $\chi$ & 0.95& 0.9 & 0.9 & & 0.95 &0.95 &0.9& &1 &0.98 &0.95 \\ [-0.5ex]
  $\gamma <$             &0.15 &0.23 &0.33 & &0.15 &0.28 &0.33& &0.2 &0.31 &0.38 \\ [-0.5ex]
\hline\hline 
\end{tabular}
\label{tab-tb-2}
\end{center}
\end{table}
A summary of of best $\chi$ result is given in Table ~\ref{tab-tb-2}. Again, a smaller $\beta$ corresponds to heavier tails and larger best $\chi$. 
Moreover, as we predicted the best $\chi$ for $\beta=3.5$, $\beta=4$ and $\beta=4.5$ 
in the range of $(0.8,1]$, $(0.67,1]$ and $(0.57,1]$, respectively.
Best $\chi$'s increase in sample size. 
\end{example}
\subsection{Combined Heavy-tailed and Long Range dependence case}
If we take $N=1$ and dimension $d=1$, we have $\displaystyle(x_k,y_{k+1}) = \sum_{j=-\infty}^\infty C_{k-j}\Xi_j,$ in which $C_j=(c_j,c_j)^T$ and 
$\Xi_j=(\xi_{j}^{(1)},\xi_{j}^{(1)})_{j \in \mathbb Z}$ are i.i.d.. 
Hence, $x_k=\sum\limits_{j\in\mathbb Z}c_{k-j}\xi_{j}^{(1)}$ and $y_{k+1}=\sum\limits_{j\in\mathbb Z}c_{k-j}\xi_{j}^{(2)}$, where $\xi_{j}^{(2)}=h\xi_{j}^{(1)}+a_j$ and 
$\{a_j\}$'s are i.i.d. zero mean random variables. 
This relation between $\xi_{j}^{(1)}$ and $\xi_{j}^{(2)}$ is due to the fact that $y_{k+1}=x_kh+\epsilon_k$ and $\epsilon_k=\sum\limits_{j\in\mathbb Z}c_{k-j}a_j$.   
We consider $\{c_j=|j|^{-\sigma}\}, \mbox{for} \,j\neq0$ and $\sigma \in (\frac12, 1 ]$, $c_0=1$. 
 The linear algorithm (\ref{a3}) reduces to:
\begin{eqnarray}\label{a-2-1}
	\quad \quad h_{k+1} = h_k +\mu_k (x_k y_{k+1}-x_k^2h_k)
	=h_k +\mu_k (x_k^2h+x_k\epsilon_k-x_k^2h_k).
\end{eqnarray}
The initial and optimal values are $h_1=401$ and $h=1$.
  
\begin{example}
Let $\xi_{j}^{(1)}\sim PL(x_{min}=0.01,\beta)$ and $a_j=f_j-E(f_j)$ with 
$f_{j} \sim PL(x'_{min}=0.01,\beta)$. The simulation is done for one-sided process and since in computer we cannot technically do infinite sum, we assume summation over the range of $(0,500,000)$. 
As in the last two examples the normalized errors in 100 trial simulations, $\{h_n^{(i)}\}_{i=1}^{100}$, are averaged 
and results for different $\chi$'s, $\beta$'s and sample sizes $n$ are presented in the following tables. The assumed  $\sigma$ is $0.65$.
The Marcinkiewicz threshold, $M=\frac{1}{\alpha}\vee(2-2\sigma)$, corresponding to $\beta=4$, $\beta=4.5$ and $\beta=5$ is
$0.7$. Hence, predicted ranges for $\chi$'s with smallest $\overline{rh}$ will be $(0.7,1]$. 
Simulation results are provided in Table ~\ref{tab:lrd-1} with summary of best $\chi$ in Table ~\ref{tab:lrd-3}.
It worth noticing that the convergence does not seem to take place below the Marcinkiewicz threshold and 
the best $\chi$s are in the predicted ranges and the normalized error decreases as $\beta$ increases.\\  

\begin{table}[h!]
\caption{Algorithm performance for LRD-HT cases with $\sigma=0.65$}
\centering  
\small\addtolength{\tabcolsep}{-4.5pt}
\begin{tabular}{l c c c c c c c c c c c }  
\hline\hline      
&\multicolumn{3}{c}{ n=100} & & \multicolumn{3}{c}{ n=5000} & &\multicolumn{3}{c}{ n=10,000} \\     
\cline{3-3} \cline{7-7} \cline{11-11}         
$\chi \backslash \beta$ & 4  & 4.5   & 5 &  & 4   & 4.5   & 5  & & 4  & 4.5   & 5 \\ 
\cline{2-4} \cline{6-8} \cline{10-12}
0.55  & \footnotesize{0.024582}  & \footnotesize{0.015503}  & \footnotesize{0.012051}  &   & \footnotesize{0.029090}  & \footnotesize{0.019016}   & \footnotesize{0.015093}    &   & \footnotesize{0.027746}  & \footnotesize{0.018037}  & \footnotesize{0.014275}     \\ 
0.6   & \footnotesize{0.010917}  & \footnotesize{0.006166}  & \footnotesize{0.004508}  &   & \footnotesize{0.013172}  & \footnotesize{0.007826}   & \footnotesize{0.005897}    &   & \footnotesize{0.012465}  & \footnotesize{0.007359}  & \footnotesize{0.005527}     \\ 
0.7   & \footnotesize{0.000665}  & \footnotesize{0.000237}  & \footnotesize{0.000132}  &   & \footnotesize{0.000958}  & \footnotesize{0.000414}   & \footnotesize{0.000262}    &   & \footnotesize{0.000881}  & \footnotesize{0.000377}  & \footnotesize{0.000238}     \\ 
0.75  & \footnotesize{2.98e-05}  & \footnotesize{7.88e-06}  & \footnotesize{6.49e-06}  &   & \footnotesize{9.77e-05}  & \footnotesize{3.15e-05}   & \footnotesize{1.70e-05}    &   & \footnotesize{8.79e-05}  & \footnotesize{2.83e-05}  & \footnotesize{1.52e-05}     \\ 
0.8   & \footnotesize{1.02e-05}  & \footnotesize{7.76e-06}  & \footnotesize{6.39e-06}  &   & \footnotesize{5.19e-06}  & \footnotesize{3.80e-06}   & \footnotesize{3.11e-06}    &   & \footnotesize{4.72e-06}  & \footnotesize{3.30e-06}  & \footnotesize{2.69e-06}     \\
0.85  & \footnotesize{9.91e-06}  & \footnotesize{7.77e-06}  & \footnotesize{6.41e-06}  &   & \footnotesize{5.01e-06}  & \footnotesize{3.91e-06}   & \footnotesize{3.21e-06}    &   & \footnotesize{4.38e-06}  & \footnotesize{3.37e-06}  & \footnotesize{2.76e-06}     \\ 
0.9   & \footnotesize{9.93e-06}  & \footnotesize{7.79e-06}  & \footnotesize{6.45e-06}  &   & \footnotesize{5.20e-06}  & \footnotesize{4.12e-06}   & \footnotesize{3.39e-06}    &   & \footnotesize{4.54e-06}  & \footnotesize{3.50e-06}  & \footnotesize{2.86e-06}     \\ 
0.95  & \footnotesize{1.00e-05}  & \footnotesize{7.90e-06}  & \footnotesize{6.55e-06}  &   & \footnotesize{5.63e-06}  & \footnotesize{4.42e-06}   & \footnotesize{3.62e-06}    &   & \footnotesize{4.73e-06}  & \footnotesize{3.65e-06}  & \footnotesize{2.98e-06}     \\ 
0.98  & \footnotesize{1.01e-05}  & \footnotesize{7.99e-06}  & \footnotesize{6.61e-06}  &   & \footnotesize{5.97e-06}  & \footnotesize{4.69e-06}   & \footnotesize{3.86e-06}    &   & \footnotesize{4.88e-06}  & \footnotesize{3.77e-06}  & \footnotesize{3.08e-06}     \\ 
1     & \footnotesize{1.02e-05}  & \footnotesize{8.04e-06}  & \footnotesize{6.65e-06}  &   & \footnotesize{6.28e-06}  & \footnotesize{4.91e-06}   & \footnotesize{4.03e-06}    &   & \footnotesize{5.02e-06}  & \footnotesize{3.89e-06}  & \footnotesize{3.19e-06}     \\ 
\hline\hline 
\end{tabular}
\label{tab:lrd-1}
\end{table}

Note that by considering $\sigma=0.65$,  the minimum of $2-2\sigma$ and $\frac1\alpha$ for all $\beta=4, 4.5$ and $5$ is $2-2\sigma$, hence we do not expect much change in the $\chi$ as $\beta$ changes. 
In addition, the rate of convergence for all considered  $\beta$'s is determined by $\gamma<\chi+2\sigma-2$.
\begin{table}[h!]
\caption{Best fixed $\chi$-Power Law, LRD with $\sigma=0.65$} 
\centering  
\begin{tabular}{c c c c c c c c c c c c}  
\hline\hline      
&\multicolumn{3}{c}{ n=100} & & \multicolumn{3}{c}{ n=5000} & &\multicolumn{3}{c}{ n=10000} \\   [-0.5ex]
\cline{3-3} \cline{7-7} \cline{11-11}     
$ \beta$ & 4  & 4.5   & 5 &  & 4   & 4.5   & 5  & & 4  & 4.5   & 5 \\   [-0.5ex]
\cline{2-4} \cline{6-8} \cline{10-12}
 \small{Best} $\chi$                   &0.85 &0.8  &0.8& &0.85 &0.8  &0.8& &0.85 & 0.8 &0.8 \\ [-0.5ex]
 \small{Resulting} $\gamma$              &0.15 &0.1 &0.1 & &0.15 &0.1 &0.1& &0.15 &0.1 &0.1 \\ [-0.5ex]
\hline\hline 
\end{tabular}
\label{tab:lrd-3}
\end{table}
\end{example} 
\section{The proof of Theorem \ref{Converge}}
{\bf Part a)} {\bf Step 1:} Reduce rate of convergence to convergence of a transformed algorithm.\\
Letting $ \eta_{k}=\left(\frac{k+1}k\right)^\gamma-1$, setting 
$g_{k}={k}^{\gamma}\left(\bar h_{k}-h\right)$ and using (\ref{b1}), one finds that
\begin{eqnarray}
g_{k+1} 
 & = &  g_{k}+\frac{1}{k^{\chi }}\left(\hat b_{k}-\hat A_{k}g_{k}\right)+\eta_kg_k,
\end{eqnarray}
where
\begin{eqnarray}
\hat b_{k} =  (k+1)^{\gamma}\left(\bar b_{k}-\bar A_{k}h\right) \mbox{  and  }
\hat A_{k} = \left(\frac{k+1}{k}\right)^\gamma\bar A_{k}.
\end{eqnarray}
However, we have by Taylor's theorem and assumption that
\begin{eqnarray}
\frac{1}{n^{\chi }}\sum\limits _{k=1}^{n}{\eta _{k}}\|A\| 
&\le& \frac{\gamma }{n^{\chi }}\sum\limits _{k=1}^{n}k^{-1} \rightarrow 0, \quad \mbox{as} \,\, n \rightarrow \infty.
\end{eqnarray}

{\bf Step 2:} Show MSLLN for new coefficients i.e. $\frac{1}{n^\chi}\sum\limits_{k=1}^{n}(\hat A_{k}-A) \rightarrow 0$, and 
$\frac{1}{n^\chi}\sum\limits_{k=1}^{n}\hat b_{k} \rightarrow 0 \quad \mbox{as} \,\, n \rightarrow \infty. $
\begin{eqnarray}\nonumber
&\!\!\!\!\!\!\!\!\!&\left\|\frac{1}{n^{\chi }}\sum\limits_{k=1}^{n}\left(\frac{k+1}{k}\right)^\gamma\left(\bar A_{k}-A\right)
-\frac{2^{\gamma }}{n^{\chi }}\sum\limits _{k=1}^{n}\left(\bar A_{k}-A\right)\right\| \\\nonumber
&\!\!\!=\!\!\!\!\!\!& \left\|\frac{1}{n^{ \chi}}\sum\limits _{k=2}^{n}\sum\limits _{j=2}^{k}\nonumber
\left[\left(\frac{j+1}{j}\right)^\gamma-\left(\frac{j}{j-1}\right)^{\gamma}\right]\left(\bar A_{k}-A\right)\right\| 
\\\nonumber
&\!\!\!\le \!\!\!\!\!\!& \sum\limits _{j=2}^{n}\left[\left(\frac{j}{j-1}\right)^{\gamma }
-\left(\frac{j+1}{j}\right)^{\gamma }\right]\frac{1}{n^{\chi }}
\left(\left\|\sum\limits _{k=2}^{n}\left(\bar A_{k}-A\right)\right\|
+\left\|\sum\limits _{k=2}^{j-1}\left(\bar A_{k}-A\right)\right\|\right)\\\nonumber
&\!\!\!\le\!\!\!\!\!& \frac{1}{n^{\chi }} \left\|\sum_{k=2}^{n}\left(\bar A_{k}-A\right)\right\| \left[2^\gamma-\left(\frac{n+1}{n}\right)^{\gamma }\right] \\\nonumber
&\!\!\!+\!\!\!\!\!\!& \sum\limits _{j=2}^{n}\left[\left(\frac{j}{j-1}\right)^{\gamma } 
-\left(\frac{j+1}{j}\right)^{\gamma }\right]\left(\frac{j-2}{n}\right)^\chi \frac{1}{(j-2)^{\chi }} \left\|\sum_{k=2}^{j-1}\left(\bar A_{k}-A\right)\right\|
\end{eqnarray}
which goes to zero by assumption and the Toeplitz lemma. By Taylor's theorem
\begin{eqnarray}\nonumber
 & &  \left|\frac{1}{n^{\chi}}\sum\limits _{k=1}^{n}(k+1)^\gamma\left(\bar b_k-\bar A_{k}h\right)
 -\frac{1}{n^{\chi}(n+1)^{-\gamma }}\sum\limits _{k=1}^{n}\left(\bar b_k-\bar A_{k}h\right)\right| \\
 & = & \frac{1}{n^{\chi}}\left|\sum\limits_{k=1}^{n-1}\sum\limits_{j=k+1}^{n}\nonumber
 \left[{j^{\gamma}}-{\left(j+1\right)^{\gamma}}\right]\left(\bar b_k-\bar A_{k}h\right)\right| \\
 & \le & \frac{1}{n^{\chi }}\sum\limits_{j=2}^{n}\gamma \nonumber
{j}^{\gamma-1 }\left|\sum\limits 
_{k=1}^{j-1}\left(\bar b_k-\bar A_{k}h\right)\right| \\
 & \le & \frac{\gamma}{n^{\chi }}\sum\limits_{j=2}^{n} 
{j}^{\chi-1 }\frac1{(j-1)^{\chi-\gamma}}\left|\sum\limits 
_{k=1}^{j-1}\left(\bar b_k-\bar A_{k}h\right)\right|,
\end{eqnarray}
which goes to zero by the Toeplitz lemma.

{\bf Step 3:} Convergence of $g_k$, hence the rate of convergence of $\bar h_k$
follows from the Proposition \ref{mainprop} with $b=0$, $\hat h_k=g_k$, $h=0$ and $\eta_k=\left(\frac{k+1}k\right)^\gamma-1$.
\eproof

\begin{proposition}\label{mainprop}
Suppose $\{\hat{A}_{k}\}_{k=1}^\infty$ is a symmetric, positive-semidefinite 
$R^{d\times d}$-valued sequence; $A$ is a (symmetric) positive-definite matrix;
$\chi \in (0,1)$; $\theta\in(\chi,1]$; $\eta_k\le \frac{\bar \eta}{k^\theta}$; $\bar \eta>0$ and 
\begin{eqnarray}\label{c1}
\hat{h}_{k+1} =\hat{h}_k + \frac{1}{k^\chi}(\hat{b}_k - \hat{A}_{k}\hat{h}_k )+\eta_k\hat{h}_k\quad\mbox{for all }\quad k=1,2,...;
\end{eqnarray}
\vspace{-0.6cm} 
\begin{eqnarray}\label{cnver}
\frac{1}{n^{\chi}} \sum_{k=1}^{n}( \hat{b}_{k}-b)\to 0 \mbox{ and } 
\frac{1}{n^{\chi}} \sum_{k=1}^{n}( \hat{A}_{k}-A)\to 0.
\end{eqnarray}
Then, $\hat{h}_n\to h\doteq A^{-1}b$ as $n\to \infty$. 
\end{proposition}
{\bf Notation:} To ease the notation in the sequel, 
we will take the product over no factors to be $1$ and the sum of no terms to be $0$. 
For convenience, we let: 
\begin{eqnarray}\label{b6}
\nu_k:=\hat{h}_k-h,\quad\quad Y_k:= \hat{A}_{k}-A, \quad\quad z_k:=\hat{b}_k-\hat{A}_{k}h.
\end{eqnarray}

\proof {\bf Step 1:} Show simplified algorithm with $A_k$'s replaced converges.\\
We note $\displaystyle\frac{1}{n^\chi} \sum_{k=1}^{n} z_{k}\to 0$ and will show 
$\nu_k\to 0$, by proving $u_k\to0$ and $w_k:=\nu_k-u_k \to 0$, where 
\begin{eqnarray}
u_{k+1}=\left(I-\frac{A}{k^\chi}+\eta_k I\right)u_k+\frac{z_k}{k^\chi}+\eta_k h\quad \mbox{subject to}\quad u_1=\nu_1.
\end{eqnarray}
By induction, we have:
\begin{equation}\label{b9}
\quad u_n=\prod_{l=1}^{n-1}\left(I-\frac{A}{l^\chi}+\eta_l I\right)u_1+\sum_{j=1}^{n-1}F_{j,n}z_j+\sum_{j=1}^{n-1}\bar{F}_{j,n}h \quad \mbox{for} \,\,n=1,2,...
\end{equation}
where
\begin{eqnarray}\label{b10}
\left\{ \begin{array}{l} F_{j,n}=\frac{1}{j^\chi} \prod_{l=j+1}^{n-1}\left(I-\frac{A}{l^\chi}+\eta_l I\right)\\
\bar{F}_{j,n}=\eta_j j^\chi F_{j,n} \,\mbox{for} j=1,2,...,n-1, n=2,3,...
\end{array} \right.
\end{eqnarray}
Hence,by (\ref{b9}), (\ref{b10}) and Lemma \ref{Two} i, ii) 
\begin{eqnarray}\label{b11}\nonumber
\lim_{n\to \infty }|u_n|&\leq&\lim_{n\to \infty }\left\|\prod_{l=1}^{n-1}\left(I-\frac{A}{l^\chi}+\eta_l I\right)\right\|
| u_1|\\
&+&\lim_{n\to \infty }\left|\sum_{j=1}^{n-1}F_{j,n}z_j\right|
+ \lim_{n\to \infty }\left|\sum_{j=1}^{n-1}\bar F_{j,n}h\right|=0.
\end{eqnarray}

{\bf Step 2:} Transfer stability from $A$ to blocks of $A_k$.\\
Define the blocks 
\begin{eqnarray}\label{b13}
\left\{ \begin{array}{l} n_k= \lfloor (a k)^ {\frac{1}{1-\chi}}\rfloor :=\max \{i\in N_0: i\leq (a k)^ {\frac{1}{1-\chi}}\} \\ I_k=\{n_k,n_k+1,\cdots,n_{k+1}-1 \}
\end{array}\right.
\end{eqnarray}
for $k=0,1,2,...$ 
and the block products
\begin{equation}\label{b16}
\quad\, U_k=\prod_{l\in I_k}\left(I-\frac{\hat{A}_l}{l^\chi}+\eta_l I\right)\,\mbox{and}\, V_{j,k}=\prod_{l=j+1}^ {n_{k+1}-1}\left(I-\frac{\hat{A}_l}{l^\chi}+\eta_l I\right)\frac{1}{j^\chi}Y_j.
\end{equation}
For the $U_k$'s we have
\small{\begin{eqnarray}\nonumber
&\!\!\!\!\!\!\!\!\!\!\!\!\!\!\!\!\!\!&\prod_{l\in I_k}\left(I-\frac{\hat{A}_l}{l^\chi}+\eta_lI\right)\!=\!
I\!-\!\sum_{l\in I_k}\frac{\hat{A}_l}{l^\chi}+\sum_{l\in I_k}\eta_lI
+\!\!\!\!\!\!\sum_{l_1,l_2\in I_k\atop{ l_1>l_2}}\!\!\!\!\left(\frac{\hat{A}_{l_1}}{l_1^\chi}-\eta_{l_1}I\right)\!\left(\frac{\hat{A}_{l_2}}{l_2^\chi}-\eta_{l_2}I\right)\\
&\!\!\!\!\!\!\!\!\!\!\!\!-\!\!\!\!\!\!\!\!\!\!\!\!&\nonumber
\sum_{l_1,l_2,l_3\in I_k\atop{l_1>l_2>l_3}}\!\!\!\!\left(\frac{\hat{A}_{l_1}}{l_1^\chi}-\eta_{l_1}I\right)\!
\left(\frac{\hat{A}_{l_2}}{l_2^\chi}-\eta_{l_2}I\right)\!\left(\frac{\hat{A}_{l_3}}{l_3^\chi}-\eta_{l_3}I\right)\!+\!\cdots(-1)^k\prod_{l\in I_k}\left(\frac{\hat{A}_l}{l^\chi}-\eta_lI\right)
\end{eqnarray}}
\normalsize
so
\begin{eqnarray}\label{ab1}\nonumber
\!\!\!\| U_k\|&\leq& \left\| I-\sum_{l\in I_k}\frac{\hat{A}_l}{l^\chi}\right\|+\sum_{l\in I_k}\eta_l
+\left\|\sum_{l_1,l_2\in I_k\atop{ l_1>l_2}}
\left(\frac{\hat{A}_{l_1}}{l_1^\chi}-\eta_{l_1}I\right)
\left(\frac{\hat{A}_{l_2}}{l_2^\chi}-\eta_{l_2}I\right)\right\|\\\nonumber
&+&\left\|\sum_{l_1,l_2,l_3\in I_k\atop{l_1>l_2>l_3}}
\left(\frac{\hat{A}_{l_1}}{l_1^\chi}-\eta_{l_1}I\right)
\left(\frac{\hat{A}_{l_2}}{l_2^\chi}-\eta_{l_2}I\right)
\left(\frac{\hat{A}_{l_3}}{l_3^\chi}-\eta_{l_3}I\right)\right\|\\
&+&\cdots+\prod_{l\in I_k}\left\|\frac{\hat{A}_l}{l^\chi}-\eta_{l}I\right\|.
\end{eqnarray}
However, we know that
$\sum_{j_1>j_2>\cdots>j_k} a_{j_1}a_{j_2}\cdots a_{j_k} \leq \frac{1}{k!}\left(\sum_{j}a_j\right)^k$ 
for $a_j\ge0$ so, it follows that
\begin{eqnarray*}
&&\!\!\!\!\!\!\!\!\!\!\!\!\!\!\!\!\!\!\sum_{l_1,l_2\in I_k\atop{ l_1>l_2}}
\left\|\frac{\hat{A}_{l_1}}{l_1^\chi}-\eta_{l_1}I\right\|
\left\|\frac{\hat{A}_{l_2}}{l_2^\chi}-\eta_{l_2}I\right\|\\\nonumber
&+&\!\!\sum_{l_1,l_2,l_3\in I_k\atop{l_1>l_2>l_3}}
\left\|\frac{\hat{A}_{l_1}}{l_1^\chi}-\eta_{l_1}I\right\|
\left\|\frac{\hat{A}_{l_2}}{l_2^\chi}-\eta_{l_2}I\right\|
\left\|\frac{\hat{A}_{l_3}}{l_3^\chi}-\eta_{l_3}I\right\|
+\cdots
+\prod_{l\in I_k}\left\|\frac{\hat{A}_l}{l^\chi}-\eta_{l}I\right\|\\\nonumber
&\leq& \sum_{m=2}^{n_{k+1}-n_k} \frac{\left(\sum\limits_{l\in I_k}\left(\frac{\| \hat{A}_l\|}{l^ \chi}+\eta_l\right)\right)^m}{m!}.
\end{eqnarray*}
As a result, we find by (\ref{ab1}) that

\begin{eqnarray}\label{b19}\nonumber
\| U_k\|&\leq& \left\| I-A\sum_{l\in I_k}\frac{1}{l^\chi}\right\|
+\left\|\sum_{l\in I_k}\frac{Y_l}{l^\chi}\right\|+\sum_{l\in I_k}\eta_l\\
&+&\sum_{m=2}^{n_{k+1}-n_k} \frac{\left(\sum\limits_{l\in I_k}\left(\frac{\| \hat{A}_l \|}{l^ \chi}+\eta_l\right)\right)^m}{m!}.
\end{eqnarray}
Now, let $\lambda_{min} \,\,\mbox{and} \,\,\lambda_{max}$ be the smallest and biggest eigenvalues of A and define $a'=\frac{a}{1-\chi}$, where $a>0$ is chosen small enough that
\begin{eqnarray}\label{b20}
a'\leq\left\{\frac{2}{\lambda_{min}+\|A\|}, \frac{1}{d\|A\|}, \frac{\lambda_{min}}{e^1(d\|A\|)^2}\right\}.
\end{eqnarray}
Then, by (\ref{b13}) and the fact that 
$$\frac{1}{1-\chi}( n_{k+1}^{1-\chi}-n_{k}^{1-\chi})\leq \sum_{l\in I_k} \frac{1}{l^\chi}\leq \frac{1}{1-\chi}( (n_{k+1}-1)^{1-\chi}-(n_{k}-1)^{1-\chi})$$
we have $\lim_{k\to \infty}\left(\sum_{l\in I_k} \frac{1}{l^\chi}-a'\right)$ is in the range of 
\begin{eqnarray*}
\left(\lim_{k\to \infty}\frac{n_{k+1}^{1-\chi}-n_{k}^{1-\chi}-a}{1-\chi}
\, \, , \, \, \lim_{k\to\infty}\frac{n_{k+1}^{1-\chi}-n_{k}^{1-\chi}-a}{1-\chi}
+\frac{n_{k}^{1-\chi}-(n_{k}-1)^{1-\chi}}{1-\chi}\right)
\end{eqnarray*} 
so by Taylor's theorem
\begin{eqnarray}\label{b21}\nonumber
\lim_{k\to\infty}\left|\sum_{l\in I_k} \frac{1}{l^\chi}-a'\right|
&\leq&\lim_{k\to \infty}\left\{\frac{1}{1-\chi}\left| n_{k+1}^{1-\chi}-
n_{k}^{1-\chi}-a\right|+\frac{1}{(n_{k}-1)^{\chi}}\right\}\\
&=& 0,
\end{eqnarray}
which also implies 
\begin{eqnarray}\label{ab2}
\lim_{k\to\infty}\sum_{l\in I_k}\eta_l\le \bar \eta \lim_{k\to\infty}n_k^{\chi-\theta} 
\sum_{l\in I_k}\frac{1}{l^\chi}= 0.
\end{eqnarray}
For arbitrary $\epsilon>0$ one finds some $K_{\epsilon}>0$ by (\ref{b21}) and (\ref{b20}) such that
\begin{eqnarray}\label{b22}\nonumber
\left\| I-A\sum_{l \in I_k}\frac{1}{l^\chi}\right\| 
&=&\max \left\{ \| A\|\sum_{l \in I_k}\frac{1}{l^\chi}-1,1-\lambda_{min}\sum_{l \in I_k}\frac{1}{l^\chi}\right\}\nonumber\\
&\leq& 1-\lambda_{min}a'+\epsilon \quad\quad\quad \mbox{for all }\quad k\geq K_{\epsilon}
\end{eqnarray}
Moreover, we can use Lemma \ref{B} of Appendix, (\ref{b6}), (\ref{cnver}), 
(\ref{b21}), (\ref{ab2}), Taylor's theorem and 
the fact $d\| A\| a'<1$ and to obtain a $K'_{\epsilon}\geq K_{\epsilon}$ such that
\begin{eqnarray}\label{b23}\nonumber
\!\!\!\!\!\!\sum_{m=2}^{n_{k+1}-n_k} \frac{\left(\sum\limits_{l\in I_k}\left(\frac{\|\hat{A}_l\|}{l^\chi}+\eta_l\right)
\right)^m}{m!} \!\!\!
&\!\!\leq\!\!\!&\!\!\!\!\!\!\sum_{m=2}^{n_{k+1}-n_k}\frac{\displaystyle(d\|A\|\sum_{l\in I_k}\frac{1}{l^\chi}+d\| \sum_{l\in I_k}\frac{Y_l}{l^\chi}\|+\sum_{l\in I_k}\eta_l)^m}{m!}\\
&\leq& e^{1+3\epsilon}\frac{(d\| A\| a'+3\epsilon)^2}{2} \quad\quad \mbox{for all}\,\, k\geq K'_{\epsilon} 
\end{eqnarray}
Therefore, by (\ref{b22}), Lemma \ref{Two} iii), (\ref{b19}) and (\ref{b23}) one finds
\begin{eqnarray}\label{b24}\nonumber
\| U_k\| \!\!\!&\!\!\!\leq\!\!\!& \left\| I-A\sum_{l\in I_k}\frac{1}{l^\chi}\right\|+\sum_{l\in I_k}\eta_l
+\left\|\sum_{l\in I_k}\frac{Y_l}{l^\chi}\right\|+\!\!\!\sum_{m=2}^{n_{k+1}-n_k} \frac{\displaystyle\left(\sum_{l\in I_k}\left(\frac{\| \hat{A}_l \|}{l^ \chi}+\eta_l\right)\right)^m}{m!}\\
&\!\!\!\leq\!\!\!& 1-\lambda_{min}a'+3\epsilon+ e^{1+3\epsilon}\frac{(d\| A\| a'+3\epsilon)^2}{2} \quad\quad \forall\,\, k\geq K'_{\epsilon}
\end{eqnarray}
Furthermore, using the fact that $ a' < \frac{\lambda_{min}}{e^1( d \| A\|)^2}$ and making for $\epsilon>0$ small enough, we find from (\ref{b24}) that, there exists a $0<\gamma<1$ and an integer $k_1>0$ such that
\begin{eqnarray}\label{b25}
\| U_k\|\leq \gamma \quad\quad \mbox{for all}\,\, k\geq k_1
\end{eqnarray}

{\bf Step 3:} Convergence of remainder $w_n$ along a subsequence using block stability of $A_k$.\\
By (\ref{c1}), (\ref{b6}), (\ref{b9}) and $w_k:=\nu_k-u_k \to 0$
\begin{eqnarray}\label{b12}
w_{n+1}=\left(I-\frac{\hat{A_n}}{n^\chi}+\eta_n I\right)w_n-\frac{1}{n^\chi}Y_n u_n \quad \mbox{for} \,\,n=1,2,\cdots
\end{eqnarray}
so it follows by (\ref{b12}) that
\begin{eqnarray}\label{b14}\nonumber
w_n&=&\prod_{l=n_k}^{n-1} \left(I-\frac{\hat{A}_l}{l^\chi}+\eta_l I\right)w_{n_k}\\
&-& \sum_{j=n_k}^{n-1}\prod_{l=j+1}^{n-1}\left(I-\frac{\hat{A}_l}{l^\chi}
+\eta_l I\right)\frac{Y_ju_j}{j^\chi} \quad \forall\quad n \ge n_k.
\end{eqnarray}
In particular,
\begin{eqnarray}\label{b15}
w_{n_{k+1}}= U_k w_{n_k}-\sum_{j\in I_k}V_{j,k}u_j\quad \mbox{for} \quad k=0,1,\cdots,
\end{eqnarray}
where $U_k$ is defined in (\ref{b16}) and 
\begin{eqnarray}\label{b16V}
V_{j,k}=\prod_{l=j+1}^ {n_{k+1}-1}\left(I-\frac{\hat{A}_l}{l^\chi}+\eta_l I\right)\frac{1}{j^\chi}Y_j.
\end{eqnarray}
By Lemma \ref{Two} v) and (\ref{b16V}) we obtain,
\begin{eqnarray}\label{b17}
\|V_{j,k}\|&\leq&\prod_{l=j+1}^{n_{k+1}-1}\left\|\left(I-\frac{\hat{A}_l}{l^\chi}+\eta_l I\right)\right\|\frac{\| Y_j\|}{j^\chi}\nonumber\\
&\leq&\prod_{l\in I_k}\left(1+\frac{\|\hat{A}_l\|}{l^\chi}+\eta_l \right)\frac{\|Y_j\|}{j^\chi}
\stackrel{j,k}\ll\frac{\| Y_j\|}{j^\chi}\quad \mbox{for} \,j\in I_k, \,\, k=0,1,\dots 
\end{eqnarray}
Therefore, by (\ref{b25}), (\ref{b17}), (\ref{b16}), (\ref{b15}), and (\ref{b3}) we have 
\begin{equation}\label{b26}
\,\,\mid w_{n_k}\mid\stackrel{k}\ll \gamma^{k-k_1}\mid w_{n_{k_1}}\mid+\sum_{l=k_1}^{k-1}\gamma^{k-l-1}\sum_{j \in I_l}\frac{\| A\|+\| \hat{A}_j\|}{j^\chi}\mid u_j\mid  \,\, \forall\, k\geq k_1.
\end{equation}
In addition,
\begin{eqnarray*}
\sum_{j \in I_l}\frac{\| A\|+\| \hat{A}_j\|}{j^\chi}\mid u_j\mid= \| A\|\sum_{j \in I_l}\frac{1}{j^\chi}\mid u_j\mid +\sum_{j \in I_l}\frac{\| \hat{A}_j\|}{j^\chi}\mid u_j\mid
\end{eqnarray*}
so using Lemma \ref{Two} iv), (\ref{b11}), (\ref{b21}) and finally applying Toeplitz Lemma, we obtain
\begin{eqnarray}\label{b27}
\lim_{l \to \infty}\sum_{j \in I_l}\frac{\| A\|+\| \hat{A}_j\|}{j^\chi}\mid u_j\mid= 0.
\end{eqnarray}
Moreover, since
\begin{eqnarray}\label{b28}
\sum_{l=k_1}^{k-1}\gamma^{k-l-1}=\frac{1-\gamma^{k-k_1}}{1-\gamma}\stackrel{k}\ll 1  \quad\quad \mbox{for all}\,\, k= k_1,k_1+1,\cdots \quad
\end{eqnarray}
it follows from (\ref{b26}), (\ref{b27}), (\ref{b28}) and the Toeplitz Lemma with 
$a_{l,k}=\gamma^{k-l-1} 1_{k_1\leq l\leq k-1}$ and $x_l=\sum_{j \in I_l}\frac{\| A\|+\| \hat{A}_j\|}{j^\chi}\mid u_j\mid$ that
\begin{eqnarray}\label{b29}\nonumber
\lim_{k \to \infty}\mid w_{n_k}\mid&\leq&\lim_{k \to \infty}\gamma^{k-k_1}\mid w_{n_{k_1}}\mid\\
 &+&\lim_{k \to \infty}\sum_{l=k_1}^{k-1}\gamma^{k-l-1}\sum_{j \in I_l}\frac{\| A\|+\| \hat{A}_j\|}{j^\chi}\mid u_j\mid=0.
\end{eqnarray}

\noindent{\bf Step 4:} Use $w_{n_k}\rightarrow 0$ to show block convergence $\max_{n\in I_k}|w_n|\rightarrow 0$.\\
Now, we return to (\ref{b14}) and find for $n\in I_k$
\begin{eqnarray} \label{b30}\nonumber
\!\!|w_n|&\!\!\!\!\!\leq\!\!\!\!\!\!& \prod_{l=n_k}^{n-1} 
\left(\left\| I-\frac{\hat{A}_l}{l^\chi}\right\|+\eta_l\right)|w_{n_k}|
+\sum_{j=n_k}^{n-1}\prod_{l=j+1}^{n-1}\!\!
\left(\left\| I-\frac{\hat{A}_l}{l^\chi}\right\|+\eta_l\right)\!\!\frac{\|Y_j\|}{j^\chi}|u_j| \\\nonumber
&\!\!\!\!\!\!\leq\!\!\!\!\!\!& \prod_{l=n_k}^{n-1}\left(1+\frac{\|\hat{A}_l\|}{l^\chi}+\eta_l\right) |w_{n_k}|
+\sum_{j=n_k}^{n-1} \prod_{l=n_k}^{n-1}\left(1+\frac{\|\hat{A}_l\|}{l^\chi}+\eta_l\right) \frac{\| Y_j\|}{j^\chi}|u_j|\\\nonumber
&\!\!\!\!\!\!\leq\!\!\!\!\!\!& \prod_{l\in I_k} \left(1+\frac{\|\hat{A}_l\|}{l^\chi}+\eta_l\right)
\left\{|w_{n_k}|+\sum_{j\in I_k}\frac{\|Y_j\|}{j^\chi}|u_j|\right\}\\
&\!\!\!\!\!\!\leq\!\!\!\!\!\!&\prod_{l\in I_k} \left(1+\frac{\| \hat{A}_l\|}{l^\chi} +\eta_l\right)
\left\{|w_{n_k}|+\sum_{j\in I_k} \frac{\|\hat{A}_j\|+\|A\|}{j^\chi}|u_j|\right\}.
\end{eqnarray}
Finally, by (\ref{b30}), (\ref{b29}), Lemma \ref{Two} v), and (\ref{b27}) we obtain
\begin{eqnarray}
\lim_{k \to\infty}\max_{n\in I_k}\mid w_n\mid=0. \quad \quad \mbox{\eproof}
\end{eqnarray} 

\noindent{\bf Part b)} By (\ref{b1}) and (\ref{b6}), $ z_k=k^\chi (\nu_{k+1}-\nu_k)+\bar{A}_k \nu_k$. 
Averaging, then reordering the sum, we have
\begin{eqnarray*}
	\frac{1}{n^\chi} \sum_{k=1}^{n} z_{k}&=&\frac{1}{n^\chi} \left(\sum_{k=1}^{n} {k^\chi}(\nu_{k+1}-\nu_k)+\sum_{k=1}^{n}\bar{A}_k \nu_k\right) {\nonumber}\\
	&=&\nu_{n+1}-\frac{1}{n^\chi} \sum_{k=1}^{n} ({k^\chi}-(k-1)^\chi)\nu_k+\frac{1}{n^\chi} \sum_{k=1}^{n} \bar{A}_k\nu_k
\end{eqnarray*}
so
\begin{eqnarray}\label{b7}\nonumber
\left|\frac{1}{n^\chi} \sum_{k=1}^{n} z_{k}\right| &\leq&
\mid \nu_{n+1}\mid+ \sum_{k=1}^{n} \frac{k^\chi-(k-1)^\chi)}{n^\chi}\mid\nu_k\mid\\
&+& \sum_{k=1}^{n} \frac{k^{\chi-1}}{n^\chi}\|\bar{A}_k\| k^{1-\chi}\mid\nu_k\mid.
\end{eqnarray}
The second and third terms on the RHS of (\ref{b7}) converge to $0$ by the 
Toeplitz lemma with $a_{n,k}=\frac{k^\chi- (k-1)^\chi}{n^\chi}$, 
$x_k=|\nu_k|$ 
and with $a_{n,k}=\frac{k^{\chi-1}\| \bar{A}_k\|}{n^\chi} $, $x_k=k^{\chi-1}|\nu_k|$ respectively. \eproof

\section{Appendix}

We first establish our promised comparison on our conditions.
\begin{lemma}\label{CondComp}
$\displaystyle\limsup_{n\to\infty}\left\|\frac{1}{n^\chi}\sum_{k=1}^n (\bar{A}_{k}- A)\right\|=0$ 
implies $\displaystyle\frac{1}{n^\chi}\sum_{k=1}^n k^{\chi-1}\|\bar{A}_{k}\|$ is bounded in $n$.
\end{lemma}
\proof
By Lemma \ref{B} (to follow) and the fact that 
$\displaystyle \sum_{k=1}^{n}k^{\chi-1}\leq \frac{n^\chi}{\chi}$, one finds that
\begin{eqnarray}\label{b4}\nonumber
\frac{1}{n^\chi} \sum_{k=1}^n k^{\chi-1}\|\bar{A}_{k}\|&\leq& \frac{d}{n^\chi}\left\|\sum_{k=1}^n k^{\chi-1}\bar{A}_{k}\right\|\\\nonumber
&\leq& \frac{d}{n^\chi}\left\|\sum_{k=1}^n k^{\chi-1}(\bar{A}_{k}-A)\right\|+\frac{d}{n^\chi}\| A\|\sum_{k=1}^n k^{\chi-1}\\
&\leq& \frac{d}{n^\chi}\left\|\sum_{k=1}^n k^{\chi-1}(\bar{A}_{k}-A)\right\|
+\frac{d\| A\|}{\chi}.
\end{eqnarray}
Hence, by the fact
$\displaystyle\sum_{j=2}^{k}(j^{\chi-1}-(j-1)^{\chi-1})=k^{\chi-1}-1$, Taylor's theorem
and the hypothesis 
\begin{eqnarray}\label{b5}\nonumber
&\!\!\!\!\!\!\!\!\!\!\!\!\!\!\!&\!\!\!\frac{1}{n^\chi} \sum_{k=1}^n k^{\chi-1}\|\bar{A}_{k}\|\\ \nonumber
&\!\!\!\!\!\!\!\!\!\leq\!\!\!&\!\!\frac{d}{n^\chi}\left\|\sum_{k=2}^n \sum_{j=2}^{k} (j^{\chi-1}-(j-1)^{\chi-1})
(\bar{A}_{k}-A)\right\|+\frac{d}{n^\chi}\left\|\sum_{k=1}^n (\bar{A}_{k}-A)\right\|+\frac{d\| A\|}{\chi}\\\nonumber
&\!\!\!\!\!\!\!\!\!\leq\!\!\!&\!\! \frac{d}{n^\chi}\left\|\sum_{j=2}^n (j^{\chi-1}-(j-1)^{\chi-1})\sum_{k=j}^{n}(\bar{A}_{k}-A)\right\|+C\\\nonumber
&\!\!\!\!\!\!\!\!\!\leq\!\!\!&\!\! d \sum_{j=2}^n \frac{(j^{1-\chi} - (j-1)^{1-\chi})}{j^{1-\chi} (j-1)^{1-\chi}}\cdot 
\frac1{n^\chi}\left\|\sum_{k=j}^{n}(\bar{A}_{k}-A)\right\|+C\\\nonumber
&\!\!\!\!\!\!\!\!\!\leq\!\!\!&\!\! d \sum_{j=2}^n \frac{(j-1)^{-\chi} (1-\chi)}{j^{1-\chi} (j-1)^{1-\chi}}. \,\, 
\frac1{n^\chi}\left(\left\| \sum_{k=1}^{n}(\bar{A}_{k}-A)\right\|
+\left\| \sum_{k=1}^{j-1}(\bar{A}_{k}-A)\right\|\right)+C\\\nonumber
&\!\!\!\!\!\!\!\!\!\leq\!\!\!&\!\! 2d(1-\chi) \sum_{j=2}^n \frac1{j^{2-\chi}}\left( 
\left\|\frac{1}{n^\chi} \sum_{k=1}^{n}(\bar{A}_{k}-A)\right\|+
\left\|\frac{1}{(j-1)^\chi} \sum_{k=1}^{j-1}(\bar{A}_{k}-A)\right\|\right)+C 
\end{eqnarray}
where $\displaystyle C=\frac{d\| A\|}{\chi}+\sup_{n}\frac{1}{n^\chi}\left\|\sum_{k=1}^n (\bar{A}_{k}-A)\right\|<\infty$.
This final term is bounded by the Toeplitz lemma and our hypothesis. 
\eproof
\\

We give our list of technical bounds used in the proof of Proposition \ref{mainprop}.
\begin{lemma}\label{Two}
Assume the setting of Proposition \ref{mainprop};
and $F_{j,k}$, $I_k$, $\{z_k\}_{k=1}^{\infty}$ and $\{Y_k\}_{k=1}^{\infty}$ are as defined 
in (\ref{b13}), (\ref{b10}) and (\ref{b6}). 
Then, following are true:
\begin{itemize}
\item[\bf i)] $\displaystyle\lim_{n\to \infty}
\left\|\prod_{l=1}^{n-1}\left(I-\frac{A}{l^\chi}+\eta_lI\right)\right\|=0$
\item[\bf ii)] $\displaystyle\lim_{n\to \infty}\left|\sum_{j=1}^{n-1}F_{j,n}z_j\right|=0$
and $\displaystyle\lim_{n\to \infty}\left|\sum_{j=1}^{n-1}\bar F_{j,n}h\right|=0$
\item[\bf iii)] $\displaystyle\lim_{k\to \infty}\left\| \sum_{l\in I_k}\frac{Y_l}{l^\chi}\right\|=0$
\item[\bf iv)] $\displaystyle\sum_{l\in I_k}\left(\frac{\|\hat{A}_l\|}{l^\chi}+\eta_l\right)\stackrel{k}\ll 1$ for all $k=0,1,\cdots$
\item[\bf v)] $\displaystyle\prod_{l\in I_k}\left(1+\frac{\| \hat{A}_l\|}{l^\chi}+\eta_l\right)\stackrel{k}\ll 1$ for all $k=0,1,\cdots$
\end{itemize}
\end{lemma}
\proof
{\bf i)} We know
\(
\left\|I-\frac{A}{l^\chi}+\eta_lI\right\| 
\)
is the maximum eigenvalue of
\(
\left((1+\eta_l)I-\frac{A}{l^\chi}\right)
\) and
\begin{eqnarray*}
0\leq \left\|\prod_{l=1}^{n-1}\left((1+\eta_l)I-\frac{A}{l^\chi}\right)\right\|\leq\prod_{l=1}^{n-1}
\left\|(1+\eta_l)I-\frac{A}{l^\chi}\right\|.
\end{eqnarray*}
Let $D\in R$ and  $\lambda_{min}>0$ be the minimum eigenvalue of $A$ and $l^*$ be large enough that: 
$1+\eta_l-\frac{\lambda_{min}}{l^\chi}>0, \, \forall\,\, l>l^*$ so 
\begin{eqnarray*}
\prod_{l=l^*}^{n-1}\left\|(1+\eta_l)I-\frac{A}{l^\chi}\right\|
&\leq&\prod_{l=l^*}^{n-1}\left(1+\frac{\bar{\eta}}{l^\theta}-\frac{\lambda_{min}}{l^\chi}\right)\\\nonumber
&\leq& \displaystyle\exp\left(\sum_{l=l^*}^{n-1}\left(\frac{\bar{\eta}}{l^\theta}
-\frac{\lambda_{min}}{l^\chi}\right)\right)\\\nonumber
&\leq& \exp\left(\int_{l^*-1}^{n-1}\frac{\bar{\eta}}{x^\theta}dx- \int_{l^*}^{n}\frac{\lambda_{min}}{x^\chi} dx\right)\\\nonumber
&\leq& \exp\left( D+ \frac{\bar{\eta}}{1-\theta} (n-1)^{1-\theta}- \frac{\lambda_{min}}{1-\chi} n^{1-\chi}\right)\\\nonumber
&\stackrel{n}\ll& 
\exp\left(\frac{-\lambda_{min}}{2-2\chi}n^{1-\chi}\right).
\end{eqnarray*}
Hence, 
\begin{eqnarray}
\prod_{l=l^*}^{n-1}\left\|(1+\eta_l)I-\frac{A}{l^\chi}\right\|\to 0 \,\, \mbox{as}\,\,n \to \infty.
\end{eqnarray}

{\bf ii)} 
\(
\|(r^\chi+\eta_rr^\chi-(r-1)^\chi)I-A\| \leq \mid (r^\chi-(r-1)^\chi)\mid
+\bar \eta r^{\chi-\theta}+\| A\|\leq 1+\bar \eta+\| A\|
\) is upper bounded $\forall r>1$ since $\chi \in (0,1)$. 
Hence, by (\ref{b10}) we have
\\
\\
\\
\\
\\
\\
\begin{eqnarray}\label{c2}\nonumber
\!\!\!\!\!\!\!\|F_{r-1,n}-F_{r,n}\|\!\!\!&\!\!\!=\!\!\!&\!\!\!\left\|\frac{1}{(r-1)^\chi} 
\prod_{l=r}^{n-1}\left((1+\eta_l)I-\frac{A}{l^\chi}\right)-\frac{1}{r^\chi} 
\prod_{l=r+1}^{n-1}\left((1+\eta_l)I-\frac{A}{l^\chi}\right)\right\|\\ \nonumber
&\!\!\!=\!\!\!&\!\!\! \left\|\prod_{l=r+1}^{n-1}\left((1+\eta_l)I-\frac{A}{l^\chi}\right)
\left[\frac{1}{(r-1)^\chi}\left((1+\eta_r)I-\frac{A}{r^\chi}\right)-\frac{1}{r^\chi}I\right]\right\|\\ \nonumber
&\!\!\!\leq\!\!\!&\!\!\! \left\| \prod_{l=r+1}^{n-1}\left((1+\eta_l)I-\frac{A}{l^\chi}\right)\right\| \frac{1}{r^\chi(r-1)^\chi}\\\nonumber
&\times&\|(r^\chi+\eta_r r^\chi-(r-1)^\chi)I-A\|\\ 
&\!\!\!\stackrel{r,n}\ll\!\!\!&\!\!\! \frac{1}{r^\chi(r-1)^\chi}\left\|\prod_{l=r+1}^{n-1}\left((1+\eta_l)I-\frac{A}{l^\chi}\right)\right\| 
\end{eqnarray}
for all $r=2,3,...,n-1,\, n=3,4,...$. Letting $\lambda$ denote an arbitrary eigenvalue of $A$
and $L^c=\left\{l:\frac{\lambda}{l^\chi}-1-\eta_l\ge c\right\}$,  
we have that
\begin{eqnarray}\label{c3}\nonumber
\!\!\!\!\!\!\left\|\prod_{l=r+1}^{n-1}\left(1+\eta_l-\frac{\lambda}{l^\chi}\right)I\right\|
&\leq& \prod_{l=r+1}^{n-1}\left(\frac{\lambda}{l^\chi}-1-\eta_l\right)
\vee\left(1+\eta_l-\frac{\lambda}{l^\chi}\right)\\\nonumber
&\leq&\prod_{l\in L^1}\left(\frac{\lambda}{l^\chi}-1-\eta_l\right)\times
\exp\left(\sum_{l=r+1,l\not\in L^0}^{n-1}\eta_l-\frac{\lambda}{l^\chi}\right)\\\nonumber
&\leq&\prod_{l\in L^1}\left(\frac{\lambda}{l^\chi}-1-\eta_l\right)
\exp\left(\sum_{l\in L^0}\frac{\lambda}{l^\chi}-\eta_l\right)
\\\nonumber
&\times& \exp\left(\sum_{l=r+1}^{n-1}\frac{\bar \eta}{l^\theta}-\frac{\lambda}{l^\chi}\right)\\\nonumber
&\stackrel{r,n}\ll& 
\exp\left(\sum_{l=r+1}^{n-1}\frac{\bar \eta}{l^\theta}-\frac{\lambda}{l^\chi}\right)\\
&\stackrel{r,n}\ll& 
\exp\left(-\frac{\lambda_{min}}{2-2\chi}\{n^{1-\chi}-(r+1)^{1-\chi}\}\right)
\end{eqnarray}
and it follows from (\ref{c3}), the fact that the eigenvectors of A span $\mathbb R^d$ and the principle of uniform boundedness that
\begin{eqnarray}\label{c4}
\left\| \prod_{l=r+1}^{n-1}\left((1+\eta_l)I-\frac{A}{l^\chi}\right)\right\|\stackrel{r,n}\ll e^{-\frac{\lambda_{min}}{2-2\chi}\{n^{1-\chi}-(r+1)^{1-\chi}\}}.
\end{eqnarray}
It follows by (\ref{b10}), (\ref{c2}) and (\ref{c4}) that
\begin{eqnarray}\label{c5}\nonumber
\sum_{r=2}^{n-1} (r-1)^{\chi}\| F_{r-1,n}-F_{r,n}\| &\stackrel{n}\ll& \sum_{r=2}^{n-1}\frac{1}{r^{\chi}}e^{-\frac{\lambda_{min}}{2-2\chi}\{n^{1-\chi}-(r+1)^{1-\chi}\}}\\\nonumber
&\stackrel{n}\ll& e^{-\frac{\lambda_{min}}{2-2\chi} n^{1-\chi}}\int_{2}^{n}\frac{1}{t^{\chi}}e^{\frac{\lambda_{min}}{2-2\chi} t^{1-\chi}}dt \\
&\stackrel{n}\ll& 1 \quad \forall\, n=3,4,...
\end{eqnarray}

Next, $\displaystyle \sum_{j=1}^{n-1}F_{j,n}z_j=\sum_{j=1}^{n-1}F_{n-1,n}z_j+\sum_{j=1}^{n-1}
\left(\sum_{r=j+1}^{n-1}F_{r-1,n}-F_{r,n}\right)z_j$  and
\begin{eqnarray}\nonumber
\left|\sum_{j=1}^{n-1}\sum_{r=j+1}^{n-1}(F_{r-1,n}-F_{r,n})z_j\right|
&\leq& \sum_{r=2}^{n-1}\| F_{r-1,n}-F_{r,n}\| \left|\sum_{j=1}^{r-1}z_j\right|.
\end{eqnarray}
Therefore, by assumption, (\ref{b10}), (\ref{c5}) and Toeplitz' lemma 
with $a_{n,r}=(r-1)^{\chi}\| F_{r-1,n}-F_{r,n}\|$ and $x_r=\frac{1}{(r-1)^\chi}\mid\sum_{j=1}^{r-1}z_j\mid$ we have: 
\begin{eqnarray}\label{c7}\nonumber
\left|\sum_{j=1}^{n-1} F_{j,n} z_j\right|&\leq& \| F_{n-1,n}\|\left|\sum_{j=1}^{n-1}z_j\right|+\left|\sum_{j=1}^{n-1}\sum_{r=j+1}^{n-1}(F_{r-1,n}-F_{r,n})z_j\right|\\
&\leq&\frac{1}{(n-1)^\chi}\left|\sum_{j=1}^{n-1}z_j\right|+\sum_{r=2}^{n-1}
\|F_{r-1,n}-F_{r,n}\| \left|\sum_{j=1}^{r-1}z_j\right|\to 0.
\end{eqnarray}
as $ n\to\infty. $Turning to the second limit in ii), we have by (\ref{b10}) and (\ref{c4}) that
\begin{eqnarray}
\sum_{j=1}^{n-1} \|\bar F_{j,n}\|\stackrel{n}\ll 
\sum_{j=1}^{n-1}j^{-\chi} 
e^{-\frac{\lambda_{min}}{2-2\chi}\{n^{1-\chi}-(j+1)^{1-\chi}\}}
j^{\chi-\theta}.
\end{eqnarray}
However, 
\begin{eqnarray}\nonumber
\sum_{j=1}^{n-1}
\frac{1}{j^{\chi}}e^{-\frac{\lambda_{min}}{2-2\chi}\{n^{1-\chi}-(j+1)^{1-\chi}\}}
&\stackrel{n}\ll& e^{-\frac{\lambda_{min}}{2-2\chi} n^{1-\chi}}
\int_{1}^{n}\frac{1}{t^{\chi}}e^{\frac{\lambda_{min}}{2-2\chi} t^{1-\chi}}dt \\
&\stackrel{n}\ll& 1 
\end{eqnarray}
for all $n$ so the second limit in ii) follows by the Toeplitz lemma.

\noindent{\bf iii)} Since
\(
\frac{1}{l^{\chi}}=\frac{1}{n_{k}^{\chi}}+\sum\limits_{r=n_k}^{l-1}\left(\frac{1}{(r+1)^{\chi}}-\frac{1}{r^{\chi}}\right)\quad \quad \forall\quad l\in I_k
\),
one has that
\begin{eqnarray*}
\left\| \sum_{l\in I_k}\frac{Y_l}{l^{\chi}}\right\|&\leq&\frac{1}{n_{k}^{\chi}}\left\|\sum_{l\in I_k}Y_l\right\|+\left\| \sum_{r<l\atop{r,l\in I_k}}\left(\frac{1}{(r+1)^{\chi}}-\frac{1}{r^{\chi}}\right)Y_l\right\|.
\end{eqnarray*}
Hence, by Taylor's theorem
\begin{eqnarray}\nonumber
\left\|\sum_{l\in I_k}\frac{Y_l}{l^{\chi}}\right\|
&\leq&\frac{1}{n_{k}^{\chi}} \left(\left\|\sum_{l<n_{k+1}}Y_l\right\| +\left\|\sum_{l<n_k}Y_l\right\|\right)\\\nonumber
&+&\sum_{r=n_k}^{n_{k+1}-2}\frac{r^\chi-(r+1)^\chi}{r^\chi (r+1)^\chi }
\left\|\sum_{l<n_{k+1}}Y_l-\sum_{l\leq r}Y_l\right\|\\\nonumber \label{c8}
&\leq&\frac{1}{n_{k}^{\chi}}\left(\left\|\sum_{l<n_{k+1}}Y_l\right\|+\left\|\sum_{l<n_{k}}Y_l\right\|\right)\\
&+&\sum_{r\in I_k}\frac{\chi}{r^{\chi+1}}\left(\left\|\sum_{l<n_{k+1}}Y_l\right\|+\left\|\sum_{l\leq r}Y_l\right\|\right),
\end{eqnarray}
where the summations all start from $l=1$ and stop at $l=n_k-1$, $r$ or $n_{k+1}-1$.
Furthermore, by the hypothesis and (\ref{b13}) we have that
\begin{eqnarray}\label{c9}
\lim_{k\to\infty}\max_{r\in I_k}\frac{1}{r^{\chi}}\left\|\sum_{l<n_{k+1}}Y_l\right\|=0
\end{eqnarray}
and the first two terms on the RHS of (\ref{c8}) go to zero.
Moreover, by (\ref{b13})
\begin{eqnarray}\label{c10}\nonumber
\sum_{r\in I_k}\frac{1}{r}&\le& \log\left(\frac{n_{k+1}-1}{n_k-1}\right)\\
&=&\log\left(\frac{\lfloor (a(k+1))^{\frac{1}{1-\chi}}\rfloor-1}{\lfloor (ak)^{\frac{1}{1-\chi}}\rfloor-1}\right)\to 0\, \mbox{as}\quad k\to\infty\quad
\end{eqnarray}
due to the fact that
 \begin{eqnarray*}
1\leq \frac{\lfloor (a(k+1))^{\frac{1}{1-\chi}}\rfloor-1}{\lfloor (ak)^{\frac{1}{1-\chi}}\rfloor-1}\leq \frac{(a(k+1))^{\frac{1}{1-\chi}}}{(ak)^{\frac{1}{1-\chi}}-2}=\frac{(\frac{k+1}{k})^{\frac{1}{1-\chi}}}{1-(\frac{2}{ak})^{\frac{1}{1-\chi}}}\to 1\, \mbox{as}\quad k\to\infty.
\end{eqnarray*}
In addition, by assumption, (\ref{c9}) and (\ref{c10})
\begin{eqnarray*}
\sum_{r\in I_k}\frac{\chi}{r^{\chi+1}}\left\|\sum_{l<n_{k+1}}Y_l\right\| \stackrel{k}\ll \sum_{r\in I_k}
\frac{1}{r} \frac{1}{n_{k}^\chi}\left\|\sum_{l<n_{k+1}}Y_l\right\|\to 0\  \mbox{as}\  k\to\infty 
\end{eqnarray*}
and
\begin{eqnarray*}
\sum_{r\in I_k}\frac{\chi}{r^{\chi+1}}\left\|\sum_{l\leq r}Y_l\right\|\leq \sum_{r\in I_k}
\frac{\chi}{r}\frac{1}{r^{\chi}}\left\|\sum_{l\leq r}Y_l\right\|\to 0 \  \mbox{as}\  k\to\infty 
\end{eqnarray*}
Hence, the last term on the RHS of (\ref{c8}) goes to zero too.\\
\noindent{\bf iv)} By Lemma \ref{B}, the fact that $\|B\|\leq\mid\mid\mid B\mid\mid\mid\leq 
\sqrt{d}\|B\|$ for a matrix with rank $d$, iii) and (\ref{b21}) we have
\begin{eqnarray}\nonumber
\sum_{l\in I_k}\left\| \frac{\hat{A}_l}{l^{\chi}}\right\| &\leq&\sum_{l\in I_k} \left|\left|\left| \frac{\hat{A}_l}{l^{\chi}}\right|\right|\right|\leq\sqrt{d}\left|\left|\left|\sum_{l\in I_k}\frac{\hat{A}_l}{l^{\chi}}\right|\right|\right|\leq d\left\|\sum_{l\in I_k} \frac{\hat{A}_l}{l^{\chi}}\right\|\\\nonumber
&\leq& d\left\|\sum_{l\in I_k} \frac{(\hat{A}_l-A)}{l^{\chi}}\right\|+ d\| A\| \sum_{l\in I_k} \frac{1}{l^{\chi}}\\
&\leq& d\left\|\sum_{l\in I_k} \frac{Y_l}{l^{\chi}}\right\|+ d\| A\| \sum_{l\in I_k} 
\frac{1}{l^{\chi}}\stackrel{k}\ll 1 \quad\mbox{for} \ k=0,1,2,...
\end{eqnarray}
Moreover, by (\ref{ab2}) $\displaystyle \sum_{l\in I_k}\eta_l\le \bar\eta \sum_{l\in I_k}\frac1{l^\theta}\to 0$
as $k\rightarrow\infty$.

\noindent{\bf v)} This follows by iv) and the fact that
\begin{eqnarray}\nonumber
\prod_{l\in I_k}\left(1+\frac{\| \hat{A}_l\|}{l^{\chi}}+\eta_l\right)&\leq& 
\exp\left(\sum_{l\in I_k}\frac{\|\hat{A}_l\|}{l^{\chi}}+\eta_l\right)\\
&\stackrel{k}\ll& 1,\, \, \quad \forall \, k=0,1,... \,\,\mbox{\eproof}
\end{eqnarray}

The following lemma is taken from Kouritzin \cite{20}.
\begin{lemma}\label{B}
Suppose $m$ is a positive integer and $\{M_k, k=1,2,3,...\}$ is a sequence of symmetric, positive 
semidefinite $R^{m \times m}$-matrices. 
Then, it follows that 
\begin{eqnarray*}
\sum_{k=1}^{j} \mid\mid\mid M_k\mid\mid\mid \leq \sqrt{m} \left|\left|\left| \sum_{k=1}^{j} M_k\right|\right|\right|, \quad\quad \forall j=1,2,3,...
\end{eqnarray*}
\end{lemma}


\begin{thebibliography}{9}
\vspace{-0.3cm}
\bibitem{12} Benveniste, A., Metivier,  M. and Priouret, P. (1990). {\it Adaptive Algorithms and Stochastic Approximation.}  New York: Springer-Verlag. 
\bibitem{31} Berbee, H.  (1987). {\it Convergence rates in the strong law for a bounded mixing sequence.} Probability Theory and Related Fields, {\bf vol}. 74, pp. 253-270.
\bibitem{13} Bertsekas, D. P. and Tsitsiklis, J. N. (1996). {\it Neuro-Dynamic Programming.} Atlanta, GA: Athena Scientific.
\bibitem{B97} Berger, E.  (1997). {\it An almost sure invariance principle for stochastic approximation procedures in linear filtering theory.} Ann. Appl. Probab. \textbf{7} pp. 444-459. 
\bibitem{27} Chen, H.-F. (1996). {\it Recent developments in stochastic approximation.} Proc. IFAC World Congr., pp. 375–380.
\bibitem{17} Chong, E. K. P., Wang, I.-J. and Kulkarni, S. R.  (1999). {\it Noise conditions for prespecified convergence rates of stochastic approximation algorithms.} IEEE Trans. Inform. Theory, {\bf vol.} 45, pp. 810–814.
\bibitem{26} Clauset, A., Shalizi, C. R. and Newman, M. E. J. (2009). {\it Power-Law Distributions in Empirical Data.} Journal SIAM Review, {\bf vol.} 51 Issue 4, pp. 661-703.
\bibitem{27} Delyon, B.  (2000). {\it Stochastic approximation with decreasing gain: Convergence and asymptotic theory.} unpublished report.
\bibitem{14} Devroye, L., Gy\"{o}rfi, L.  and Lugosi, G. (1996). {\it A Probabilistic Theory of Pattern Recognition.} Berlin, Germany: Springer-Verlag.
\bibitem{22} Dippon \"{J}.  and Walk, H.  (2006). {\it The Averaged Robbins – Monro Method for Linear Problems in a Banach Space.} Journal of Theoretical Probability, {\bf vol.} 19, No. 1.
\bibitem{18} Eweda E. and Macchi, O. (1984). {\it Convergence of an adaptive linear estimation algorithm.} IEEE Trans. Automat. Contr., {\bf vol.} AC-29, pp. 119–127.
\bibitem{23} Even-Dar, E. and Mansour, Y.  (2004). {\it Learning rates for q-learning.} Journal of Machine Learning Research, {\bf vol.} 5, pp. 1–25.
\bibitem{5}  Farden, D. C.  (1981). {\it Stochastic Approximation with Correlated Data.} IEEE Trans. Inform. Theory, {\bf vol.} IT-27, NO. 1.
 \bibitem{8} Frost, O. L.  (1972). {\it An algorithm for linearly constrained adaptive array processing.} Proc. IEEE, {\bf vol.} 60, pp. 922-935.
 \bibitem{24} George, A. P.  and Powell, W. B. (2006). {\it Adaptive stepsizes for recursive estimation with applications in approximate dynamic programming.} Journal of Machine Learning Research, {\bf vol.} 65, pp. 167–198.
 \bibitem{7} Griffiths, L. J. (1969). {\it A simple algorithm for real-time processing in antenna arrays.} Proc. IEEE, {\bf vol.} 57, pp. 1696-1704.
 \bibitem{2} Gy\"{o}rfi,  L.  (1980). {\it Stochastic approximation from ergodic sample for linear regression.} Z. Wahrscheinlichkeitstheorie und verwandte Gebiete, {\bf vol.} 54, pp. 47–55.
 \bibitem{19} Gy\"{o}rfi, L.  (1984). {\it Adaptive linear procedures under general conditions.} IEEE Trans. Inform. Theory, {\bf vol.} IT-30, pp. 262-267.
 \bibitem{28} Karagiannis, T. , Molle, M.  and Faloutsos,  M.  (2004). {\it Long-Range Dependence Ten Years of Internet Traffic Modeling.} IEEE Computer Society. 
 \bibitem{MS-small} Kouritzin, M. A.  and Sadeghi, S.  {\it Marcinkiewicz Law of Large Numbers for Covariances of Heavy-tailed, Long-range Dependent Data.} submitted.
 \bibitem{6} Kouritzin,  M.A. (1994). {\it Inductive methods and rates of r-mean convergence in adaptive filtering.} Stochastics and Stochastics Reports 51, pp. 241-266.
 \bibitem{20}  Kouritzin, M. A.  (1996). {\it On the convergence of linear stochastic approximation procedures.} IEEE Trans. Inform. Theory, {\bf vol.} 42, pp. 1305–1309.
 \bibitem{MK96}  Kouritzin,  M. A.  (1996). {\it On the interrelation of almost sure invariance principles for  certain stochastic adaptive algorithms and for partial sums of random variables.} J. Theoret. Probab. Theory, {\textbf 42} 811–840.
 \bibitem{3} Kushner, H. J.  and Yin,  G. (2003). {\it Stochastic Approximation and Recursive Algorithms and Applications.} Springer, Second edition, pp. 8.
 \bibitem{15} Ljung, L., Pflug, G.  and  Walk, H.  (1992). {\it Stochastic Approximation and Optimization of Random Systems.} Basel, Switzerland: Birkhäuser-Verlag.
 \bibitem{30} Louhchi, S.  and Soulier, P.  (2000). {\it Marcinkiewicz-Zegmond Strong Laws for Infinite Variance Time Series.} Statistical Inference for Stochastic Processes, {\bf vol}. 3, pp. 31-40.
 \bibitem{29} Rio, E.  (1995). {\it A Maximal Inequality and Dependent Marcinkiewicz-Zegmond Strong Laws.} The Annals of Probability, {\bf vol.}23 No. 2, pp. 918-937.
 \bibitem{4} Robbins, H. and Monro,  S.  (1951). {\it A stochastic approximation method.} Ann. Math. statist., {\bf vol.} 22, pp. 400-407.
 \bibitem{16} Solo,  V.  and  Kong, X.  (1995). {\it Adaptive Signal Processing Algorithms: Stability and Performance.} Englewood Cliffs, NJ: Prentice-Hall.
 \bibitem{25} Stout, W. F.  (1974). {\it Almost Sure Convergence.} Academic Press Inc., pp. 126.
 \bibitem{11} \ Tadi\'{c}, V.B. (2004). {\it On the Almost Sure Rate of Convergence of Linear Stochastic Approximation Algorithms.} IEEE Trans. Inform. Theory, {\bf vol}. 50, No. 2.
 \bibitem{21} Walk, H.  and Zsid\'{o},  L. (1989). {\it Convergence of Robbins-Monro method for linear problems in banach space.} J. Math. Anal. Applic., {\bf vol}. 139, pp. 152–177.
 \bibitem{1} Yin, G.  (1992). {\it Asymptotic Optimal Rate of Convergence for an Adaptive Estimation Procedure}  Stochastic Theory and Adaptive Control.
 
\end{thebibliography}
\end{document}